%% file: MOS_paper.tex
\documentclass[journal,twocolumn]{IEEEtran}

\usepackage{amsthm}
\usepackage{amsmath}
\usepackage{amsfonts}
\usepackage{dsfont}
\usepackage{mathrsfs}
\usepackage{pifont}
\usepackage{amssymb}
\usepackage{verbatim}
\usepackage{upgreek}
\usepackage{color}
\usepackage[final]{graphicx}
\usepackage{caption}
\usepackage{algorithmic}
\usepackage{algorithm}
\usepackage{epsfig}
\usepackage{eucal}
\usepackage{cite}
\usepackage{subfigure}

\makeatletter
\def\maketag@@@#1{\hbox{\m@th\normalfont\normalsize#1}}
\makeatother

\input{Defs2}

\title{Model Order Selection Rules For Covariance Structure Classification}

\vspace{0.1cm}

\author{V. Carotenuto, \IEEEmembership{Member, IEEE}, A. De Maio, \IEEEmembership{Fellow, IEEE}, D. Orlando, \IEEEmembership{Senior Member, IEEE}, P. Stoica, \IEEEmembership{Fellow, IEEE}
\thanks{V. Carotenuto and A. De Maio are with the Dipartimento di Ingegneria Elettrica e delle Tecnologie dell'Informazione, Universit\`a degli Studi di Napoli ``Federico II'', via Claudio 21, I-80125 Napoli, Italy. E-mail: {\tt ademaio@unina.it}, {\tt vincenzo.carotenuto@unina.it}.}
\thanks{D. Orlando is with Universit\`a degli Studi ``Niccol\`o Cusano'', via Don Carlo Gnocchi 3, 00166 Roma, Italy. 
E-mail: {\tt danilo.orlando@unicusano.it}.}
\thanks{P. Stoica is with the Department of Information Technology, Uppsala University, P O Box 337, SE-751 05, Uppsala, Sweden. 
E-mail: {\tt peter.stoica@it.uu.se}.}
}

\begin{document}

\maketitle

\begin{abstract}
The adaptive classification of the interference covariance matrix structure for radar signal processing applications is addressed in this paper. This represents a key issue because many detection architectures are synthesized assuming a specific covariance structure which may not necessarily coincide with the actual 
one due to the joint action of the system and environment uncertainties.  The considered classification problem is cast in terms of a multiple hypotheses test with some nested alternatives and the theory of Model Order Selection (MOS) is exploited to devise suitable decision rules. 
Several MOS techniques, such as the Akaike, Takeuchi, and Bayesian information criteria are adopted and the corresponding merits and drawbacks are discussed.
At the analysis stage, illustrating examples for the probability of correct model selection are presented showing the effectiveness of the proposed rules.
\end{abstract}

\section{Notation}
In the sequel, vectors and matrices are denoted by boldface lower-case and upper-case letters, respectively.
The symbols $\det(\cdot)$, $\tr(\cdot)$, $\otimes$, $(\cdot)^*$, $(\cdot)^T$, $(\cdot)^\dag$ denote the determinant, trace, Kronecker product, complex conjugate, transpose, and conjugate transpose, respectively. As to numerical sets, $\R$ is the set of real numbers, $\R^{N\times M}$ is the Euclidean space of $(N\times M)$-dimensional real matrices (or vectors if $M=1$), $\C$ is the set of complex numbers, and $\C^{N\times M}$ is the Euclidean space of $(N\times M)$-dimensional complex matrices (or vectors if $M=1$). The symbols $\Re\left\{ z \right\}$ and $\Im\left\{ z \right\}$ indicate the real and imaginary parts of the complex number $z$, respectively. $\bI_N$ stands for the $N \times N$ identity matrix, while $\bzero$ is the null vector or matrix of proper dimensions. We denote by $\bJ\in\R^{N\times N}$ a permutation matrix such that $\bJ(l,k)=1$ if and only if $l+k=N+1$. Given a matrix $\bA=[\boa_1,\ldots,\boa_M] \in \C^{N \times M}$, $\vect(\bA) = [\boa_1^T, \boa_2^T, \ldots, \boa_M^T]^T \in \C^{NM \times 1}$, while given a vector $\boa\in\C^{N\times 1}$, $\diag(\boa)\in\C^{N\times N}$ indicates 
the diagonal matrix whose $i$th diagonal element is the $i$th entry of $\boa$.

The Euclidean norm of a vector is denoted by $\|\cdot\|$. 
We write $\bM \succ \bzero$ if $\bM$ is positive definite. 
Let $f(\bx)\in\R$ be a scalar-valued function of vector argument, then $\partial f(\bx)/\partial \bx$ denotes the gradient of $f(\cdot)$ with respect to $\bx$ arranged
in a column vector, while $\partial f(\bx)/\partial \bx^T$ is its transpose. Moreover, if $\widehat{\bx}$ belongs to the domain of $f(\cdot)$, then the gradient
of $f(\cdot)$ with respect to $\bx$ and evaluated at $\widehat{\bx}$ is denoted by $\partial f(\widehat{\bx})/\partial \bx$. 
For a finite set $A,\; |A|$ stands for its cardinality. $U(N)\subset \C^{N\times N}$ denotes the set of all $N \times N$ unitary matrices and $j=\sqrt{-1}$. 
For two sets, $A$ and $B$, $A \times B$ denotes their Cartesian product. The $(k,l)$-entry (or $l$-entry) of a 
generic matrix $\bA$ (or vector $\boa$) is denoted by $\bA(k,l)$ (or $\boa(l)$). 
Given two statistical hypotheses $H_i$ and $H_j$, then $H_i\subset H_j$ means that $H_i$ is nested into $H_j$.
The acronym i.i.d. means independent and identically 
distributed while the symbol $E[\cdot]$ denotes statistical expectation.
Finally, we write $\bx\sim\cC\cN_N(\bm, \bM)$ if $\bx$ is a complex circular $N$-dimensional normal vector with mean $\bm$ and covariance matrix $\bM \succ \bzero$, $\bx\sim \cN_N(\bm, \bM)$ if $\bx$ is a $N$-dimensional normal vector with mean $\bm$ and covariance matrix $\bM \succ \bzero$, and 
$\varphi\sim \cU(0,2\pi)$ if $\varphi$ is a random variable uniformly distributed in $(0,2\pi)$.

\section{Introduction, Motivation, and Problem Formulation}
\label{Sec:ProbForm}

Consider a radar system equipped with $N\geq 2$ (spatial and/or temporal) channels. The echoes from the cell under test (CUT) are downconverted
to baseband, pre-processed, properly sampled, and organized to form a $N$-dimensional vector, $\bz$ say, referred to as primary data or CUT sample. 
A set of secondary data, $\bz_1,\ldots,\bz_{K}$, with $K > N$, statistically independent of $\bz$, is also acquired 
in order to make the system adaptive with respect to the unknown
Interference Covariance Matrix (ICM), $\bM\succ\bzero$.
As is customary, these data are assumed to share the same ICM as $\bz$ and are obtained exploiting 
echoes from range cells in the proximity of the CUT within the reference window \cite{kelly1986adaptive,robey1992cfar,BOR-Morgan,Cai1992,
DeMaioSymmetric,Pascal,JunLiu00,JunLiu01,JunLiu02,CP00,DeMaio-RAO}.

To accomplish the detection task which is typical of the search process, the radar signal processor 
solves a testing problem applying 
a decision rule computed from the collected data (decision statistic).
From a mathematical viewpoint, target detection can be formulated in terms of a binary hypothesis test and tools
provided by the {\em Decision Theory} can be exploited to solve it. Several design criteria have been adopted in this respect: 
the Generalized Likelihood Ratio Test (GLRT) \cite{kelly1986adaptive,Yuri01,GLRT-based,DD,Raghavan}, 
the Wald test \cite{HaoPHE,DeMaio-WALD-AMF,WLiuRao,GCuiRaoWald,GCuiRaoWaldCompound,WLiuRao}, 
the Rao test\footnote{Note that GLRT, Wald test, and Rao test, under mild conditions, are asymptotically equivalent \cite{KayBook}.}
\cite{DeMaio-RAO,GCuiRaoWald,GCuiRaoWaldCompound,JunLiu00,WLiuRao,RicciRao}, and 
the Invariance Principle \cite{BoseSteinhardtInvariance,DeMaioInvCoinc,DeMaioInvPersymmetry,RaghavanInv,Ciuonzo1,Ciuonzo2}.

Usually a given design technique is applied under specific assumptions on the ICM structure which are tantamount to incorporating some degree of {\em a priori} knowledge at the design stage.
Specifically, certain structures of the covariance $\bM$ can be induced by
the interference type, the geometry of the system array, and/or uniformity of the transmitted pulse train.
In the most general case, $\bM\in\C^{N\times N}$ is Hermitian, but it is well-known that:
\begin{itemize}
\item ground clutter, observed by a stationary monostatic radar, often exhibits a symmetric power spectral density centered around the 
zero-Doppler frequency implying that the resulting ICM is real, i.e., $\bM\in\R^{N \times N}$ \cite{Klemm};
\item from a theoretical point of view, symmetrically spaced linear arrays or pulse trains induce a persymmetric structure on 
$\bM$ \cite{Nitzberg-1980}; the following two cases are possible
\begin{itemize}
\item $\bM\in\C^{N\times N}$ is Hermitian and persymmetric (or centrohermitian) if and only if $\bM=\bJ\bM^*\bJ$;
\item $\bM\in\R^{N\times N}$ is symmetric and persymmetric (or centrosymmetric) if and only if $\bM=\bJ\bM\bJ$.
\end{itemize}
\end{itemize}

\begin{figure}[htp!]
\begin{center}
\includegraphics[width=\columnwidth]{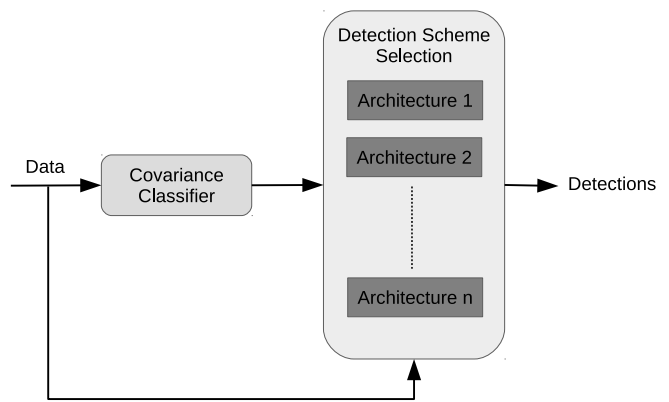}
\caption{Block diagram of a two-stage detection architecture exploiting the covariance structure classifier.}
\label{fig:figure01}
\end{center}
\end{figure}

For each of the mentioned scenarios, there exist examples of 
adaptive detectors in the literature \cite{DeMaioSymmetric,Cai1992,HaoSP_HE}. The
knowledge about the environment as well as the structure of the ICM can guide 
the system operator towards the most appropriate
decision scheme. In this regard, the primary sources of available information are directly related to the system and/or to the operating scenario.
However, there exist a plethora of causes that introduce uncertainty and make the nominal assumptions no longer valid.
For instance, array calibration errors would produce residual imbalances among channels that can heavily degrade the ICM persymmetric structure. 
Another example concerns the level of symmetry of ground clutter power spectral density which can be altered by
the possible presence of a dominating Doppler or some discretes with a given velocity.
This motivates the need for a classifier capable of inferring the ICM structure over the range bins of the system reference window.
Its output could then be fed to a selector choosing the most suitable detection scheme as shown in Figure \ref{fig:figure01}.

A possible approach to handle the mentioned classification problem is based on its formulation in terms of a multiple hypothesis test and on the use of model order selection (MOS) rules, since each 
possible choice for $\bM$ represents a {\em model} with a given number of 
parameters \cite{StoicaBabu1,Stoica1,StoicaSelenLi,Kay1,Kay2,Kay3,Kay2005,kay2005multifamily}. 
Following this idea, it is worth making explicit the relationship between parameters and model. 
To this end, note that the number of parameters  
introduced by 
the specific structure of $\bM$ can be stacked into a vector $\btheta_{i} \in \mathbb{R}^{m_{i}\times 1}$, 
where $m_{i}$ depends on the specific scenario. 
Since the entries of $\btheta_{i}$ parameterize $\bM$,
this dependence is denoted using the notation $\bM(\btheta_{i})$.
Finally, the considered models (or hypotheses) are representative of combinations among the possible assumptions on 
the clutter spectrum (symmetry around zero-Doppler or the lack of the mentioned symmetry) and
the system configuration (persymmetry).

In summary, the problem at hand is tantamount to choosing among the following 
hypotheses:
\begin{equation} \label{multihyptest}
\begin{cases}
H_1: \bM(\btheta_{1})\in\C^{N\times N} &\mbox{is Hermitian unstructured,}\\
H_2: \bM(\btheta_{2})\in\R^{N\times N} &\mbox{is symmetric unstructured,}\\
H_3: \bM(\btheta_{3})\in\C^{N\times N} &\mbox{is centrohermitian,}\\
H_4: \bM(\btheta_{4})\in\R^{N\times N} &\mbox{is centrosymmetric.}\\
\end{cases}
\end{equation}
The number of unknown parameters under each hypothesis is given by:
\be 
\begin{cases}
m_1 = N^2         & \mbox{under } H_1,\\
m_2 =N(N+1)/2    & \mbox{under } H_2,\\
m_3 =N(N+1)/2    & \mbox{under } H_3,\\

m_4 =
\begin{cases}
\dfrac{N}{2} \left(\dfrac{N}{2}+1\right)               & \mbox{if } N \mbox{ is even}\\
\left( \dfrac{N+1}{2} \right)^2 & \mbox{if } N \mbox{ is odd}\\
\end{cases}
& \mbox{under } H_4.
\end{cases}
\label{eqn:ParameterNumber}
\ee

For the sake of clarity, the proofs of \eqref{eqn:ParameterNumber} for the cases centrohermitian and centrosymmetric 
are provided in Appendix \ref{App:PersymmetricPar}.

Hereafter, for brevity, we omit the dependence on $\btheta_i$ letting 
$\bM_i = \bM(\btheta_i)$ and $\bX_i = \bM^{-1}(\btheta_i)$.

Before concluding this section, a few remarks are in order. First, notice that different models could have the same number of parameters but, as 
shown in the next sections, this is not a limitation since classification rules exploit specific estimates 
corresponding to the different structures reflecting the assumed hypothesis. 
Second, it is possible to identify nested hypotheses among those listed in \eqref{multihyptest}, for instance $H_2\subset H_1$, $H_3\subset H_1$,
$H_4\subset H_2$, etc.

In the next section, several MOS classification algorithms for problem \eqref{multihyptest} are briefly described highlighting the respective design
assumptions, which might not be always met in the considered radar application. 
The latter observation means that the behavior of these classification rules 
versus the parameters of interest deserves a careful investigation. Section \ref{Sec:Derivation} provides closed-form expressions for the classification
statistics discussed in Section \ref{Sec:MOS_description}. Concretely, these statistics are computed according to two approaches.
The first exploits the overall data matrix which also comprises the CUT, whereas the second neglects the CUT and uses secondary data only.
The performances
of the considered selectors are analyzed in Section \ref{Sec:Performance}, where the figure of merit is 
the probability of correct classification as a function of the number of data used for estimation. 
Finally, concluding remarks and future research tracks are given in Section \ref{Sec:Conclusions}.
Mathematical derivations are confined to the appendices.

\section{Model Order Selection Criteria}
\label{Sec:MOS_description}
The aim of this section is twofold. The first part provides useful preliminary definitions, while the second part
presents a brief review of the adopted selection criteria for problem \eqref{multihyptest}. Subsequent developments assume that
$\bz_{k}\sim\cC\cN_N(\bzero,\bM)$, $k=1,\ldots,K$, and $\bz\sim\cC\cN_N(\alpha \bv,\bM)$, where $\alpha=\alpha_{re}+j\alpha_{im}$, $\alpha_{re}, \alpha_{im}\in\R$, is 
an amplitude factor accounting for target response and propagation effects and 
$\bv\in\C^{N\times 1}$ is the nominal steering vector. Finally, the vectors $\bz_1,\ldots,\bz_K,\bz$ are assumed to be statistically independent. 

Now, denote by $\bZ = \left[\bz_1,\ldots,\bz_K \right] \in \mathbb{C}^{N \times K}$ the entire secondary data matrix and let
$\bp_i$ be the parameter vector under the $H_i$ hypothesis, $i=1,\ldots,4$. Observe that
\begin{itemize}
\item if the CUT is incorporated into the classification rules, then $\bp_i=[\btheta_{i}^T  \ \balpha^T]^T \in\R^{n_i \times 1}$, where
$\balpha=[\alpha_{re} \ \alpha_{im}]^T\in\R^{2\times 1}$, $n_i=m_i+2$; in this case, we let $Z_c=\{ \bz, \bZ \}$;
\item if the the classification rules are devised from $\bZ$ only, then $\bp_i=\btheta_{i} \in\R^{n_i}$, where $n_i=m_i$; here we let $Z_c=\{ \bZ \}$.
\end{itemize}
Because the derivation of the MOS criteria requires the computation of the maximum likelihood estimates (MLE) of the unknown parameters as well as suitable estimates of the Fisher Information Matrix (FIM), let us provide the expressions of the probability density functions (pdfs) of $\bz$, $\bz_k$, $k=1,\ldots,K$, $\bZ$, and the joint pdf of $\bz$ and $\bZ= \left[\bz_1,\ldots,\bz_K \right] \in \mathbb{C}^{N \times K}$ under the considered hypotheses, namely, $\forall i=1,\ldots,4$:
\begin{equation} \label{eqn:pdfz}
f(\bz;\bp_i, H_i) = \frac{\exp\left\{-(\bz-\alpha\bv)^\dag\bX_i(\bz-\alpha\bv)\right\}}{\pi^N \det(\bM_i)},
\end{equation}
\begin{equation} \label{eqn:pdfzk}
f(\bz_k;\btheta_i, H_i) = \frac{\exp\left\{-\bz_k^\dag\bX_i\bz_k\right\}}{\pi^N \det(\bM_i)}, \quad k=1,\ldots,K,
\end{equation}
\begin{equation}
\begin{split}
f \left(\bZ ;\bp_i, H_i \right) & = \prod_{k=1}^K f(\bz_k;\btheta_i, H_i)\\
& = 
\left\{ \dfrac{\exp \left\lbrace -\frac{1}{K} \tr \left[ \bX_i  \bS \right] \right\rbrace}
{\pi^N \det(\bM_i)} \right\}^{K}
\end{split}
\end{equation}
\begin{equation} \label{eqn:pdfzZ}
\begin{split}
f (\bz, \bZ ; & \bp_i, H_i ) = f(\bz;\bp_i, H_i) \prod_{k=1}^K f(\bz_k;\btheta_i, H_i)\\
&=
\left\{ \dfrac{\exp \left\lbrace -\frac{1}{K+1} \tr \left[ \bX_i ( \bS_{\alpha} + \bS) \right] \right\rbrace}{\pi^N \det(\bM_i)} \right\}^{K+1}
\end{split}
\end{equation}

where $\bS_{\alpha}=(\bz-\alpha\bv)(\bz-\alpha\bv)^\dag$ and $\bS=\bZ \bZ^\dag$.

Finally, denote by $s(\bp_{i}, H_i ; \bz)=\log f \left( \bz;\bp_i, H_i \right)$, $s(\btheta_{i}, H_i ; \bz_k)=\log f \left( \bz_k; \btheta_i, H_i \right)$, $k=1,\ldots,K$, and let
\be
s(\bp_{i},  H_i ; Z_c)=
\begin{cases}
s(\bp_{i},  H_i ; \bz, \bZ) = \log f \left( \bz, \bZ ;\bp_i, H_i \right),\\ 
\mbox{if the CUT is included},\\
s(\bp_{i},  H_i ; \bZ) = \log f \left( \bZ ;\bp_i, H_i \right),\\
\mbox{if the CUT is excluded},\\
\end{cases}
\ee
denote the log-likelihood functions\footnote{Observe that $\alpha$ is a nuisance parameter with respect to problem (\ref{multihyptest}).}.

The remainder of this section is focused on MOS criteria. Several of such criteria have been developed for the selection of 
an estimated best approximating model from a set of candidates \cite{Anderson};
most of them rely on minimization of the Kullback-Leibler (KL) discrepancy. A well-known rule is the Akaike Information Criterion (AIC), which, with reference to
problem \eqref{multihyptest}, can be formulated as
\be
H_{\widehat{i}}=
\arg\min_{\H} \left\{
-2s(\widehat{\bp}_{i}, H_i ; Z_c)+2n_i, 
\right\}, \quad\quad\quad \mbox{(AIC)}
\label{eqn:AIC}
\ee
where $H_{\widehat{i}}$ is the estimated model, $\H=\{ H_1,\ldots,H_4 \}$, and $\widehat{\bp}_{i}$ is the MLE of $\bp_i$. 
The main drawback of this rule is its non-zero probability of overfitting \cite{Stoica1} due to the penalty term $2n_i$ being too small
for high-order models, especially for nested hypotheses. 
To overcome this limitation, an empirical modification of AIC has been proposed in \cite{BhansaliGIC}. This rule, referred to as
Generalized Information Criterion (GIC), corrects the penalty term of AIC via a factor $(1+\rho)$ with $\rho>1$, namely
\be
H_{\widehat{i}}=
\arg\min_{\H} \left\{
-2s(\widehat{\bp}_{i}, H_i ; Z_c)+(1+\rho)n_i 
\right\} \quad\quad\quad \mbox{(GIC)}.
\label{eqn:GIC}
\ee
Note that if we set $\rho=1$ GIC reduces to AIC. 

The Takeuchi Information Criterion (TIC), whose main goal is to extend AIC to mismodeling scenarios, has the following form \cite{Anderson}:
\be
H_{\widehat{i}}=\arg\min_{\H} \left\{
-2s(\widehat{\bp}_{i}, H_i ; Z_c)+
2\tr[ \widehat{\cbJ}_i(\widehat{\bp}_{i})  \widehat{\cbI}_i^{-1}(\widehat{\bp}_{i})]
\right\}, \quad\quad\quad \mbox{(TIC)}
\label{eqn:TIC}
\ee
where $\cbI_i(\widehat{\bp}_{i})\in\R^{n_i\times n_i}$ is the negative Hessian of the log-likelihood function 
evaluated at $\widehat{\bp}_{i}$, namely the observed FIM, whose expression is
\be
\widehat{\cbI}_i(\widehat{\bp}_{i})=-\frac{\partial^2 s(\widehat{\bp}_{i}, H_i ; Z_c)}{\partial \bp_i \partial \bp_i^T},
\ee
and $\widehat{\cbJ}_i(\widehat{\bp}_{i})$ is the sample FIM, viz.
\begin{equation}
\begin{split}
\widehat{\cbJ}_i(\widehat{\bp}_{i})=&
\frac{\partial s(\widehat{\bp}_{i}, H_i ; \bz)}{\partial \bp_i}
\frac{\partial s(\widehat{\bp}_{i}, H_i ; \bz)}{\partial \bp^T_i}\\
& +
\sum_{k=1}^K 
\left[
\frac{\partial s(\widehat{\btheta}_{i}, H_i ; \bz_k)}{\partial \bp_i}
\frac{\partial s(\widehat{\btheta}_{i}, H_i ; \bz_k)}{\partial \bp^T_i}
\right].
\end{split}
\end{equation}
when $\bz$ and $\bZ$ are both considered or
\be
\sum_{k=1}^K 
\left[
\frac{\partial s(\widehat{\btheta}_{i}, H_i ; \bz_k)}{\partial \bp_i}
\frac{\partial s(\widehat{\btheta}_{i}, H_i ; \bz_k)}{\partial \bp^T_i}
\right].
\ee
when only $\bZ$ is considered. Note that, given the true model parameter vector $\bar{\bp}$ and the true hypothesis $\bar{H}$, 
$\widehat{\cbI}_i(\widehat{\bp}_{i})$ and $\widehat{\cbJ}_i(\widehat{\bp}_{i})$ are estimators of
\be
{\cbI}(\bar{\bp})=-E\left[\frac{\partial^2 s(\bar{\bp}, \bar{H} ;Z_c)}{\partial \bar{\bp} \partial \bar{\bp}^T}\right]
\ee
and
\be
{\cbJ}(\bar{\bp})=
E\left[
\frac{\partial s(\bar{\bp}, \bar{H} ;Z_c)}{\partial \bar{\bp}}
\frac{\partial s(\bar{\bp}, \bar{H} ;Z_c)}{\partial \bar{\bp}^T}
\right],
\ee
respectively. It is important to observe that, in general, ${\cbI}(\bar{\bp})$ will not equal ${\cbJ}(\bar{\bp})$ when the model is misspecified.
However, if the model is correctly specified, then by the {\em Information Matrix Equivalence Theorem} \cite{AIC_Bozdogan}, the information matrix
can be expressed in either Hessian form, ${\cbI}(\bar{\bp})$, or in the outer product form, ${\cbJ}(\bar{\bp})$. 

Both the AIC (along with its generalization) and TIC are derived under the assumption of large samples. To relax
this requirement, the corrected AIC (AICc) has been devised:
\be
H_{\widehat{i}}=\arg\min_{\H} \left\{
-2s(\widehat{\bp}_{i}, H_i ;Z_c)+
2n_i \frac{(K+1)N}{(K+1)N-n_i-1}
\right\} \quad\quad\quad \mbox{(AICc)}.
\label{eqn:AICc}
\ee
It is important to note that in the considered framework the AICc is essentially a heuristic rule since it
has been originally proposed for linear regression models \cite{Sugiura} and later extended to 
the case of nonlinear regression and autoregressive time series \cite{hurvich}, which neither covers the scenarios considered herein.

Finally, other selection rules, such as the Bayesian Information Criterion (BIC), can be obtained according to a Bayesian framework.
The BIC has been derived as an asymptotic approximation to a transformation of the Bayesian posterior probability 
of a candidate model \cite{Schwarz}.
In large-sample settings, BIC selects the model which is {\em a posteriori} most probable.
It is also worth mentioning that, under some regularity conditions, BIC minimizes the KL discrepancy \cite{Anderson,Stoica1}.
An alternative formulation of BIC can be obtained relaxing the large-sample requirement and assuming a noninformative prior for both the
parameter vector $\btheta_i$ and the model $H_i$. Under the above hypotheses, BIC can be expressed as \cite{Stoica1,StoicaBabu,BICneath}
\be
H_{\widehat{i}}=\arg\min_{\H} \left\{
-2s(\widehat{\bp}_{i}, H_i ; Z_c)+\log[\det(\widehat{\cbI}_i(\widehat{\bp}_{i}))]
\right\}, \quad\quad\quad \mbox{(BIC)}
\label{eqn:BIC}
\ee
which, for large samples and the herein considered context, reduces to (see Subsection \ref{subs:AsymptBIC})
\be
H_{\widehat{i}}=\arg\min_{\H} \left\{
-2s(\widehat{\bp}_{i}, H_i ; Z_c)+m_i\log(K)
\right\}. \quad\quad\quad \mbox{(Asymptotic BIC)}
\label{eqn:BIC_asympt}
\ee
We note, once again, that even though different models can share the same number of parameters, the considered selection criteria are still capable of discriminating between the different hypotheses since they use the specific MLEs together with the corresponding log-likelihood function under the current hypothesis.

Also note that the definition of large or small samples, which is important for some of the previous criteria, depends 
on the ratio between the number of parameters, $n_i$, and number 
of data, $(K+1)N$ or $KN$. Moreover, for the considered application, $n_i$ depends on $N$.
Thus, the behavior of these criteria might change according to the specific application and, for this reason, has to be investigated.

For the problem under consideration, the ratio between the number of parameters and the number of samples approaches zero as the number of 
homogeneous secondary data, $K$, increases. However, this situation might not be realizable in practical scenarios with the consequence that
the large samples assumption would be no longer valid. Finally, the presence of outliers, clutter-edges, and/or regions with
highly varying reflectivity can make the assumption that the true model belongs to the family of candidates fail. 
Thus, given these
uncertainty factors, it is worthwhile investigating the considered MOS rules to determine
which one performs better than the others. This is the scope of the next sections.

\section{Computation of MOS Decision Rules}
\label{Sec:Derivation}
This section contains the derivation of the explicit expressions of the aforementioned classification rules.
Specifically, we follow two approaches: Approach A jointly exploits secondary and primary data; whereas Approach B relies on secondary data only.
The the former processes an additional data vector (primary data) with respect to the latter, but the number of unknown
parameters increases due to the presence of the target complex amplitude. Moreover, the estimate of the target response represents an additional
computational load for the rules based on the full data, which requires
the computation of the decision statistics for each look direction. In contrast to this, 
Approach B does not depend on the 
system steering vector and, hence, the classification schemes can be evaluated
irrespective of the current steering direction. The above strategies are described in the next two subsections, whereas
the last subsection provides the expression of BIC for large values of $K$.

\subsection{MOS Decision Rules Using the Entire Data Matrix}
It follows from Section \ref{Sec:MOS_description} that the ingredients needed to construct a MOS decision rule are the MLEs of the unknown parameters, 
the log-likelihood functions, and the matrices
$\widehat{\cbI}_i(\bp_i)$ and $\widehat{\cbJ}_i(\bp_i)$. 
Evidently the mathematical expressions for all the above quantities depend on which model ($H_i$) is assumed.

The log-likelihood functions can be easily obtained from \eqref{eqn:pdfz}, \eqref{eqn:pdfzk}, and \eqref{eqn:pdfzZ}, namely
\begin{equation} \label{eqn:scorez}
s(\bp_{i}, H_i ;\bz) = - N \log \pi - \log \det(\bM_i) - \tr\left\lbrace \bX_i\bS_{\alpha}\right\rbrace,
\end{equation}
\begin{equation} \label{eqn:scorez_k}
s(\btheta_{i}, H_i ;\bz_k) = - N\log \pi - \log \det(\bM_i) - \tr \left\lbrace \bX_i \bS_{k} \right\rbrace \, , 
\end{equation}
\begin{equation} \label{eqn:scorezZ}
\begin{split}
s(\bp_{i}, H_i ;\bz, \bZ) = & - (K+1) \left[ N \log \pi + \log \det(\bM_i) \right]\\ 
& - \tr \left\lbrace \bX_i \bS \right\rbrace -\tr\left\lbrace \bX_i \bS_\alpha \right\rbrace,
\end{split}
\end{equation}
where $\bS_k=\bz_k\bz_k^\dag$.

The next step towards the derivation of the MOS statistics consists 
in evaluating the gradients of $s(\bp, H_i; \bz)$ and $s(\bp, H_i; \bz_k)$, $k=1,\ldots,K$, which are required to compute $\widehat{\cbJ}_i(\bp_i)$. 
More precisely, observe that
\be
\frac{\partial s(\bp_{i}, H_i ; \bz)}
{\partial \bp_i} =
\left[
\begin{array}{c}
\ds\frac{\partial s(\bp_{i}, H_i ; \bz)}
{\partial \btheta_i}
\\
\vspace{-5mm}
\\
\ds\frac{\partial s(\bp_{i}, H_i ; \bz)}
{\partial \balpha}
\end{array}
\right]
\ee
and
\be
\frac{\partial s(\btheta_{i}, H_i ; \bz_k)}
{\partial \bp_i} =
\left[
\begin{array}{c}
\ds\frac{\partial s(\btheta_{i}, H_i ; \bz_k)}
{\partial \btheta_i}
\\
\ds\bzero
\end{array}
\right].
\ee
In Appendix \ref{App:Grad_z_zk}, it is shown that
\begin{equation} \label{eqn:gradient_z}
\begin{split}
& \frac{\partial s(\bp_{i}, H_i ; \bz)}{\partial \btheta_i} =\\
&
\left\{
\begin{array}{l}
-\{\{\vect[\bX_i]\}^\dag \bC_i\}^T + \bC_i^\dag \left[\bX_i^* \otimes \bX_i \right] \vect[\bS_\alpha], 
\\
\mbox{if } \bM_i \mbox{ is Hermitian},
\\
-\{\{\vect[\bX_i]\}^T \bC_i\}^T + \bC_i^T (\bX_i \otimes \bX_i)\vect[\bS_\alpha], 
\\
\mbox{if } \bM_i \mbox{ is symmetric},
\end{array}
\right.
\end{split}
\end{equation}
where $\bC_i \in \C^{N^2 \times m_i}$ is a transformation matrix that depends on the specific structure of 
$\bM_i$ and on how $\btheta_{i}$ is defined (see also Appendix \ref{App:Grad_z_zk}),
\begin{equation} \label{eqn:gradient_zk}
\begin{split}
& \frac{\partial s(\btheta_{i}, H_i ; \bz_k)}{\partial \btheta_i} =\\
& \left\{
\begin{array}{l}
-\{\{\vect[\bX_i]\}^\dag \bC_i\}^T + \bC_i^\dag \left[\bX_i^* \otimes \bX_i \right] \vect[\bS_k], 
\\
\mbox{if } \bM_i \mbox{ is Hermitian},
\\
-\{\{\vect[\bX_i]\}^T \bC_i\}^T + \bC_i^T (\bX_i \otimes \bX_i)\vect[\bS_k], 
\\
\mbox{if } \bM_i \mbox{ is symmetric},
\end{array}
\right.
\end{split}
\end{equation}
and
\begin{equation} \label{eqn:gradient_zAalpha}
\frac{\partial s(\bp_{i}, H_i ; \bz)}{\partial \balpha}
= 2 
\begin{bmatrix}
-\alpha_{re} \bv^\dag \bX_i \bv +\Re\left\{ \bz^\dag \bX_i\bv \right\}\\
-\alpha_{im} \bv^\dag \bX_i \bv -\Im\left\{ \bz^\dag \bX_i\bv \right\}\\
\end{bmatrix} .
\end{equation}
Now, we move to the evaluation of the Hessian of $s(\bp_i,H_i; \bz, \bZ)$, which can be partitioned as follows
\begin{align}
& \widehat{\cbI}_i(\bp_i) =
- \frac{\partial^2 s(\bp_{i}, H_i ; \bz, \bZ)}{\partial \bp_{i}\bp_{i}^T} \nonumber
\\
& = -
\begin{bmatrix}
\ds
\frac{\partial^2 s(\bp_{i}, H_i ;\bz, \bZ)}{\partial \btheta_{i}\btheta_{i}^T} & 
\ds
\frac{\partial^2 s(\bp_{i}, H_i ;\bz, \bZ)}{\partial \btheta_{i}\balpha^T}
\\
\vspace{-4mm}
\\
\ds
\frac{\partial^2 s(\bp_{i}, H_i ;\bz, \bZ)}{\partial \balpha\btheta_{i}^T} & 
\ds
\frac{\partial^2 s(\bp_{i}, H_i ;\bz, \bZ)}{\partial \balpha\balpha^T}
\end{bmatrix}\nonumber
\\
&=
-
\begin{bmatrix}
\bH_{\theta\theta,i} & \bH_{\alpha\theta,i}^T
\\
\bH_{\alpha\theta,i} & \bH_{\alpha\alpha,i}
\end{bmatrix},
\end{align}
where
\begin{equation} \label{eqn:HessThetaTheta}
\begin{split}
& \bH_{\theta\theta,i}=\\
& \begin{cases}
\bC_i^\dag \{\bX_i^* \otimes [(K+1)\bX_i - \bX_i(\bS+\bS_\alpha)\bX_i]\\
-[\bX_i(\bS+\bS_\alpha)\bX_i]^* \otimes \bX_i\}\bC_i,
\mbox{ if } \bM_i \mbox{ is Hermitian},
\\
\vspace{-4mm}
\\
\bC_i^T \{\bX_i \otimes [(K+1)\bX_i - \bX_i(\bS+\bS_\alpha)\bX_i]\\
- \bX_i(\bS+\bS_\alpha)^* \bX_i \otimes \bX_i\}\bC_i, \mbox{ if } \bM_i \mbox{ is symmetric},
\end{cases}
\end{split}
\end{equation}
$\bH_{\alpha\alpha,i}=-2\bv^{\dag} \bX_i\bv \bI_2$, and if $\bM_i$ is Hermitian
\be
\bH_{\alpha\theta,i}=\left[\begin{array}{c}
\left\{
2\alpha_{re} \bC_i^\dag [\bX_i^*\otimes \bX_i]\vect[\bv\bv^\dag]\right.
\\
\left.-2\Re\{ \bC_i^\dag [\bX_i^*\otimes \bX_i] \vect[\bv\bz^\dag] \}
\right\}^T
\\
\vspace{-4mm}
\\
\left\{
2\alpha_{im} \bC_i^\dag [\bX_i^*\otimes \bX_i]\vect[\bv\bv^\dag]\right.
\\
\left. +2\Im\{ \bC_i^\dag [\bX_i^*\otimes \bX_i] \vect[\bv\bz^\dag] \}
\right\}^T
\end{array}\right]
\label{eqn:HessAlphaThetaHerm}
\ee
while if $\bM_i$ is symmetric
\be
\bH_{\alpha\theta,i}=\left[\begin{array}{c}
\left\{
2\alpha_{re} \bC_i^T [\bX_i \otimes \bX_i]\vect[\bv\bv^\dag]\right.
\\
\left. -2\Re\{ \bC_i^T [\bX_i \otimes \bX_i] \vect[\bv\bz^\dag] \}
\right\}^T
\\
\vspace{-4mm}
\\
\left\{
2\alpha_{im} \bC_i^T [\bX_i \otimes \bX_i]\vect[\bv\bv^\dag]\right.
\\
\left. +2\Im\{ \bC_i^T [\bX_i \otimes \bX_i] \vect[\bv\bz^\dag] \}
\right\}^T
\end{array}\right].
\label{eqn:HessAlphaThetaSymm}
\ee
The proofs of the above statements are provided in Appendix \ref{App:Hess}.

The final step consists in replacing the unknown parameters, namely $\alpha$ and $\btheta_i$, with suitable 
estimates. Forasmuch as the ML estimates of the unknown parameters are not always available in closed form 
(to our best knowledge), we replace them with consistent estimates as follows. For the ICM, we use the ML 
estimates obtained from secondary data only. As to $\alpha$, its estimate is obtained according to the ML rule assuming 
known ICM and, then, replacing the ICM with the corresponding consistent estimate.
Thus, when the ICM is unstructured, namely under
$H_1$, the estimates of the $\bM$ and $\alpha$ are \cite{robey1992cfar}
\be
\widehat{\bM}_1 = \dfrac{1}{K} \bZ \bZ^\dag, \quad
\widehat{\alpha} = \dfrac{\bv^\dag \widehat{\bM}_1^{-1}\bz}{\bv^\dag \widehat{\bM}_1^{-1}\bv},
\label{eqn:sampleH1}
\ee
respectively.
When $H_2$ is assumed, the ICM is unstructured and real. Thus, following the lead of \cite{DeMaioSymmetric}, we use the following estimates
\be
\widehat{\bM}_2 = \dfrac{1}{K}
\Re \left\lbrace \bZ\bZ^\dag \right\rbrace  \,
\ee
and
\begin{equation}
\widehat{\alpha}_{re} =
\frac{  \Re\{\bv\}^T \widehat{\bM}_2^{-1} \Re\{\bz\} + \Im\{\bv\}^T \widehat{\bM}_2^{-1} \Im\{\bz\} }
{  \Re\{\bv\}^T \widehat{\bM}_2^{-1} \Re\{\bv\} + \Im\{\bv\}^T \widehat{\bM}_2^{-1} \Im\{\bv\}  },
\end{equation}
\begin{equation}
\widehat{\alpha}_{im} =
\frac{   \Re\{\bv\}^T \widehat{\bM}_2^{-1} \Im\{\bz\} - \Im\{\bv\}^T \widehat{\bM}_2^{-1} \Re\{\bz\}   }
{   \Re\{\bv\}^T \widehat{\bM}_2^{-1} \Re\{\bv\} + \Im\{\bv\}^T \widehat{\bM}_2^{-1} \Im\{\bv\}   }.
\end{equation}
The persymmetric structure of the ICM, which occurs under $H_3$, yields the following estimates \cite{Cai1992}
\be
\widehat{\bM}_3 = \dfrac{1}{2K} \left[ 
\bZ \bZ^\dag + \bJ (\bZ \bZ^\dag)^{*} \bJ \right]\, ,
\ee
\be
\widehat{\alpha}_{re} =
\frac{  \bv^\dag \widehat{\bM}_3^{-1} \bz_e }
{  \bv^\dag \widehat{\bM}_3^{-1} \bv }, \quad \mbox{and} \quad
\widehat{\alpha}_{im} =-j\frac{  \bv^\dag \widehat{\bM}_3^{-1} \bz_o }
{  \bv^\dag \widehat{\bM}_3^{-1} \bv },
\ee
where $\bz_e=(\bz+\bJ\bz^*)/2$ and $\bz_o=(\bz-\bJ\bz^*)/2$. 

Finally, the estimates under $H_4$ can be obtained exploiting the results in \cite{HaoSP_HE}, namely
\be
\widehat{\bM}_4 = \dfrac{1}{2K} \Re \left\lbrace
\bZ \bZ^\dag + \bJ (\bZ \bZ^\dag)^{*} \bJ \right\rbrace \, ,
\label{eqn:sampleH4}
\ee
\begin{equation}
\widehat{\alpha}_{re} =
\frac{  \tr[\bV^\dag \widehat{\bM}_4^{-1} \bZ_e ] }
{ \tr[ \bV^\dag \widehat{\bM}_4^{-1} \bV ] }, 
\quad
\widehat{\alpha}_{im} =-j\frac{  \tr[\bV^\dag \widehat{\bM}_4^{-1} \bZ_o ]}
{ \tr[ \bV^\dag \widehat{\bM}_4^{-1} \bV ]},
\end{equation}
where $\bV=[\Re\{\bv\} \ \Im\{\bv\}]$, $\bZ_e=[\Re\{\bz_e\} \ \Im\{\bz_e\}]$, and $\bZ_o=[\Re\{\bz_o\} \ \Im\{\bz_o\}]$.

\subsection{MOS Decision Rules Using Secondary Data Only}

Here we derive the expressions for the terms needed to compute the MOS rules based on secondary data only. To this end, we rely on
the previous results. More precisely, first recall that $\bp_i=\btheta_i$, $i=1,\ldots,4$, and the log-likelihood function is given by (see \eqref{eqn:scorez_k}) 
\begin{equation}
\begin{split}
& s(\btheta_{i}, H_i ; \bZ) =\\
& - K \left[ N \log \pi + \log \det(\bM_i) \right] - \tr \left\lbrace \bX_i \bS \right\rbrace.
\end{split}
\end{equation}
Moreover, the observed FIM and the sample FIM become
\be
\widehat{\cbI}_i(\widehat{\btheta}_{i})=-\frac{\partial^2 s(\widehat{\btheta}_{i}, H_i ; \bZ)}{\partial \btheta_i \partial \btheta_i^T},
\ee
and
\be
\widehat{\cbJ}_i(\widehat{\btheta}_{i})=
\sum_{k=1}^K 
\left[
\frac{\partial s(\widehat{\btheta}_{i}, H_i ; \bz_k)}{\partial \btheta_i}
\frac{\partial s(\widehat{\btheta}_{i}, H_i ; \bz_k)}{\partial \btheta^T_i}
\right].
\ee
respectively, where $\widehat{\btheta}_{i}$ is the ML estimate of $\btheta_i$ under $H_i$. 
Note that, as opposed to Approach A, in this case closed form expressions for the 
ML estimates are available and they are precisely given by the expressions presented in the previous subsections (see \eqref{eqn:sampleH1}-\eqref{eqn:sampleH4}).
Finally, to evaluate the gradient of $s(\bp_{i}, H_i ; \bz_k)$, we can use \eqref{eqn:gradient_zk} and for
the Hessian of $s(\btheta_{i}, H_i ; \bZ)$, we use \eqref{eqn:HessThetaTheta} after replacing $\bS+\bS_\alpha$ with $\bS$.

\subsection{BIC for Large K}\label{subs:AsymptBIC}
In this subsection, we specialize \eqref{eqn:BIC} in the limit of $K\rightarrow +\infty$. To this end, we first consider Approach A and
approximate the penalty term of BIC as
\begin{align}
& \log \det[\widehat{\cbI}_i(\bp_i)] = \nonumber \\
& \log\det[ -\bH_{\theta\theta,i} ]+\log\det[-\bH_{\alpha\alpha,i} 
+\bH_{\alpha\theta,i} \bH^{-1}_{\theta\theta,i}\bH_{\theta\alpha,i} ] \nonumber
\\
&=m_i\log(K+1) + \log\det \left[ -\frac{\bH_{\theta\theta,i}}{(K+1)} \right] \nonumber\\
& +\log\det\left[-\bH_{\alpha\alpha,i} 
+\bH_{\alpha\theta,i} \frac{(K+1)}{K+1} \bH^{-1}_{\theta\theta,i}\bH_{\theta\alpha,i} \right] \nonumber
\\
&\stackrel{K\rightarrow +\infty}{\approx} m_i\log(K) + \cO(1),
\label{eqn:penaltyBIC}
\end{align}
where $\cO(1)$ represents a term that tends to a constant as $K\rightarrow +\infty$. 
The limiting approximation in \eqref{eqn:penaltyBIC} was obtained using the following asymptotic equalities
\begin{equation}
\frac{1}{K+1}(\bS+\bz\bz^\dag) \stackrel{K\rightarrow +\infty}{\approx} \bM,
\quad
\frac{1}{K}\bS \stackrel{K\rightarrow +\infty}{\approx} \bM,
\end{equation}
in the expression of $\bH_{\theta\theta,i}/(K+1)$, see \eqref{eqn:HessThetaTheta}, and observing 
that $\bH_{\alpha\alpha,i}$, \eqref{eqn:HessAlphaThetaHerm}, and \eqref{eqn:HessAlphaThetaSymm} do not depend on $K$.
As a consequence, the following equalities hold
\be
\lim_{K\rightarrow +\infty} \frac{\bH_{\alpha\theta,i}}{K+1}=\lim_{K\rightarrow +\infty} \frac{\bH_{\theta\alpha,i}}{K+1}=0
\ee
\begin{align}
\lim_{K\rightarrow +\infty} \frac{-\bH_{\theta\theta,i}}{K+1} = \bC,
\end{align}
where $\bC\succ \bzero$ does not depend on $K$. Therefore, neglecting the $\cO(1)$ term,
\eqref{eqn:BIC} becomes \eqref{eqn:BIC_asympt}.
Observe that the above criterion is also valid in the case where the CUT is not used (i.e., Approach B). As a matter of fact, the expression
of asymptotic BIC for the latter case can be obtained considering $\bH_{\theta\theta,i}$ only and repeating the above arguments.

\section{Numerical Examples and Discussion}
\label{Sec:Performance}
This section is devoted to the analysis of the classification schemes presented in the previous sections. The metric used to assess
their performance is the Probability of Correct Classification ($P_{cc}$) estimated under each hypothesis by means of standard 
Monte Carlo counting techniques over $1000$ independent trials.

The interference is modeled as circular complex normal random vectors with the following covariance matrix
\be
\bM_i = \bA_i \bR_i \bA_i^\dag + \sigma^2_n\bI, \quad i=1,\ldots,4
\ee
where $\sigma^2_n\bI$ represents the thermal noise component with $\sigma^2_n$ being its power, $\bR_i$ accounts for
the clutter contributions and incorporates the clutter power, and $\bA_i$ is a matrix factor modeling possible array channel errors as, for instance,
amplification and/or delay errors, calibration residuals, and mutual coupling \cite{Klemm}. 
The specific instances of $\bA_i$ and $\bR_i$ depend on which hypothesis is in force as shown below.

Different interference sources (with exponentially shaped covariance) are encompassed by $\bR_i$, whose $(h,k)$th entry has the following expression
\be
\bR(h,k) = {\ds \sum_{i=1}^L} \mbox{CNR}_l \rho_l^{|h-k|} e^{j 2 \pi (h-k) f_l },
\ee
where, for the $l$th interference source, CNR$_l>0$ is the Clutter-to-Noise Ratio, $\rho_l$ is the one-lag correlation coefficient, 
and $f_l$ is the normalized Doppler frequency. Finally, $L$ is the number of interference sources. For each hypothesis, 
we choose $\bR_i$ and $\bA_i$ as follows
\begin{itemize}
\item under $H_1$: $\bA_1 = \bI + \sigma_d \bW_1$, $f_l \neq 0$, $\forall l=1,\ldots,L$, where $\sigma_d>0$ and $\bW_1(h,k)\sim\cC\cN_1(0, 1)$ i.i.d.;
\item under $H_2$: $\bA_2 = \bI + \sigma_d \bW_2$, $f_l = 0$, $\forall l=1,\ldots,L$, where $\sigma_d>0$ and $\bW_2(h,k)\sim\cN_1(0, 1)$ i.i.d.;
\item under $H_3$: $\bA_3 = \bI$, $f_l \neq 0$, $\forall l=1,\ldots,L$;
\item under $H_4$: $\bA_4 = \bI$, $f_l = 0$, $\forall l=1,\ldots,L$.
\end{itemize}
As to the target signature, we choose $\alpha = \sqrt{\mbox{SNR}} e^{j \varphi}$ with $\varphi\sim \cU(0,2\pi)$ 
and SNR$ = 10$ dB is the Signal-to-Noise Ratio, whereas, the steering vector $\bv$ is chosen such 
that 
\be
\bv=\frac{1}{\sqrt{N}}\left[e^{-j2\pi f_v\frac{(N-1)}{2}} \ \cdots \ e^{-j2\pi f_v} 1 \ e^{j2\pi f_v} \ \cdots \ e^{j2\pi f_v\frac{(N-1)}{2}}\right]^T
\ee
assuming $N$ odd and $f_v=0.01$. Finally, two study cases are considered:
Case 1 assumes $L=1$, i.e., only one clutter source is considered; Case 2 considers $L=2$, i.e., two clutter types with different powers are assumed.
The latter case can arise in scenarios where the radar swath contains an edge separating two types of clutter sources (e.g., ground and sea clutter).
The considered parameter settings are described in Table \ref{tab:case_study}. 

\begin{figure}[H] \centering
\subfigure[Hypothesis 1.]{\includegraphics[width=0.49\columnwidth]{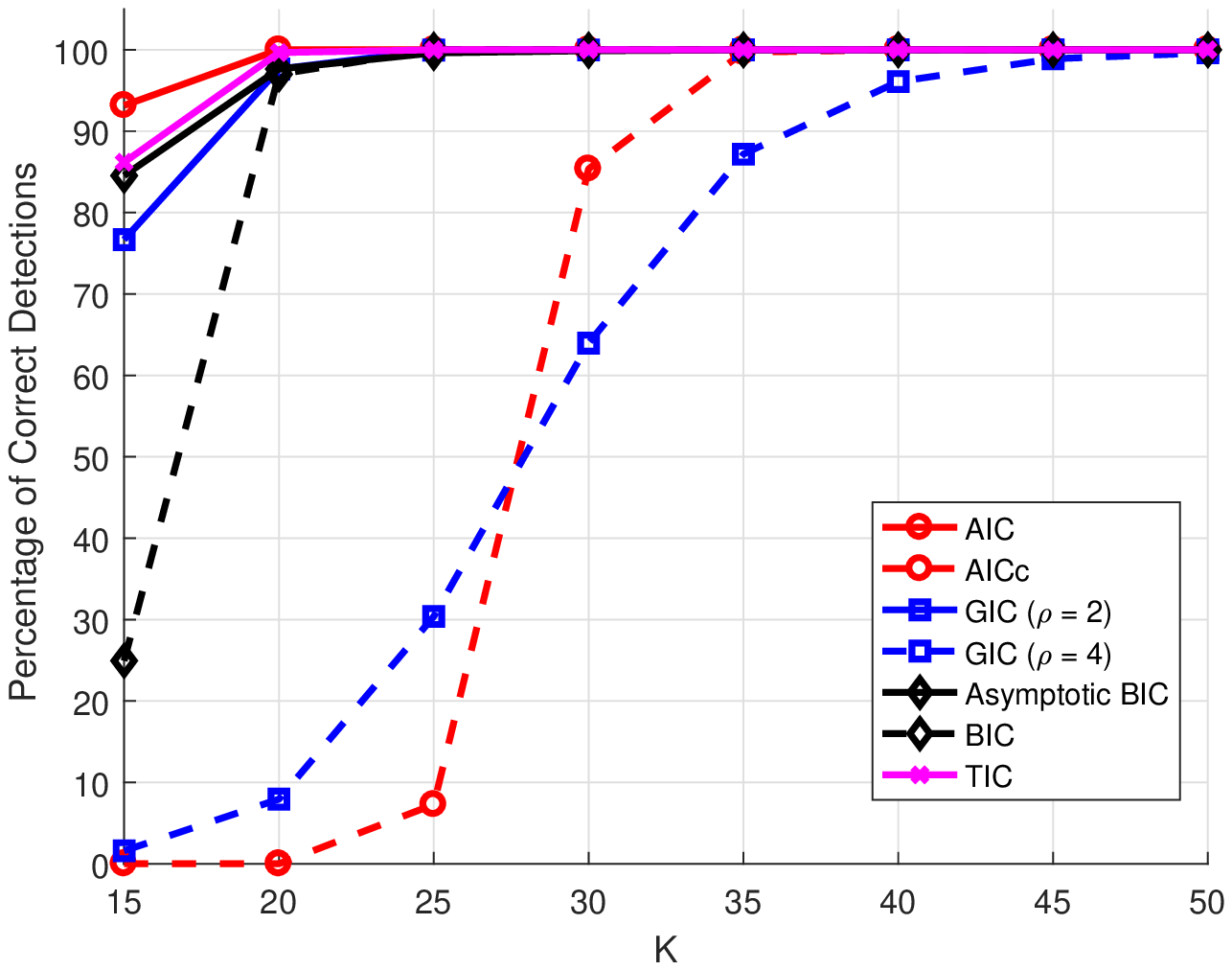}}
\subfigure[Hypothesis 2.]{\includegraphics[width=0.49\columnwidth]{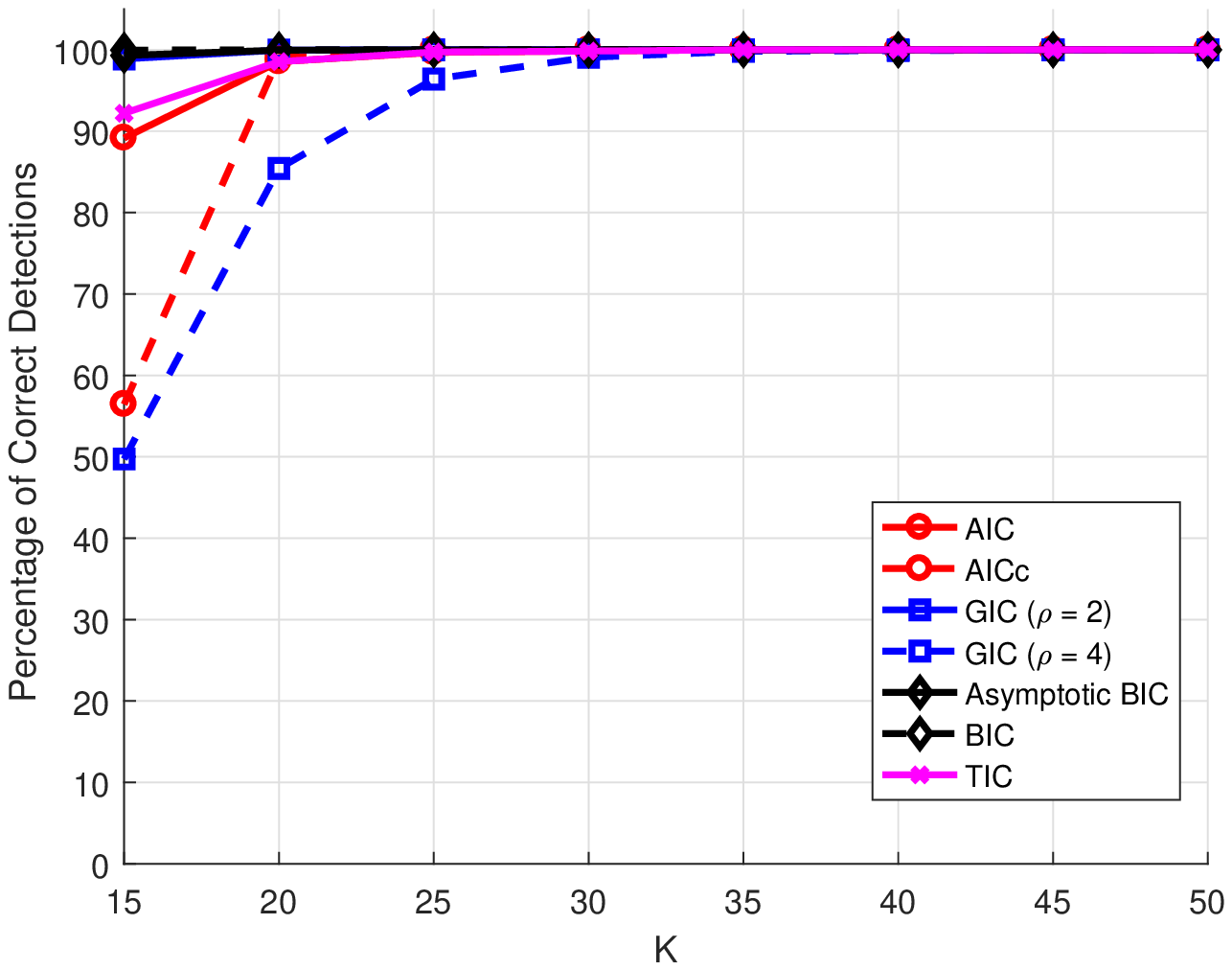}}\\
\subfigure[Hypothesis 3.]{\includegraphics[width=0.49\columnwidth]{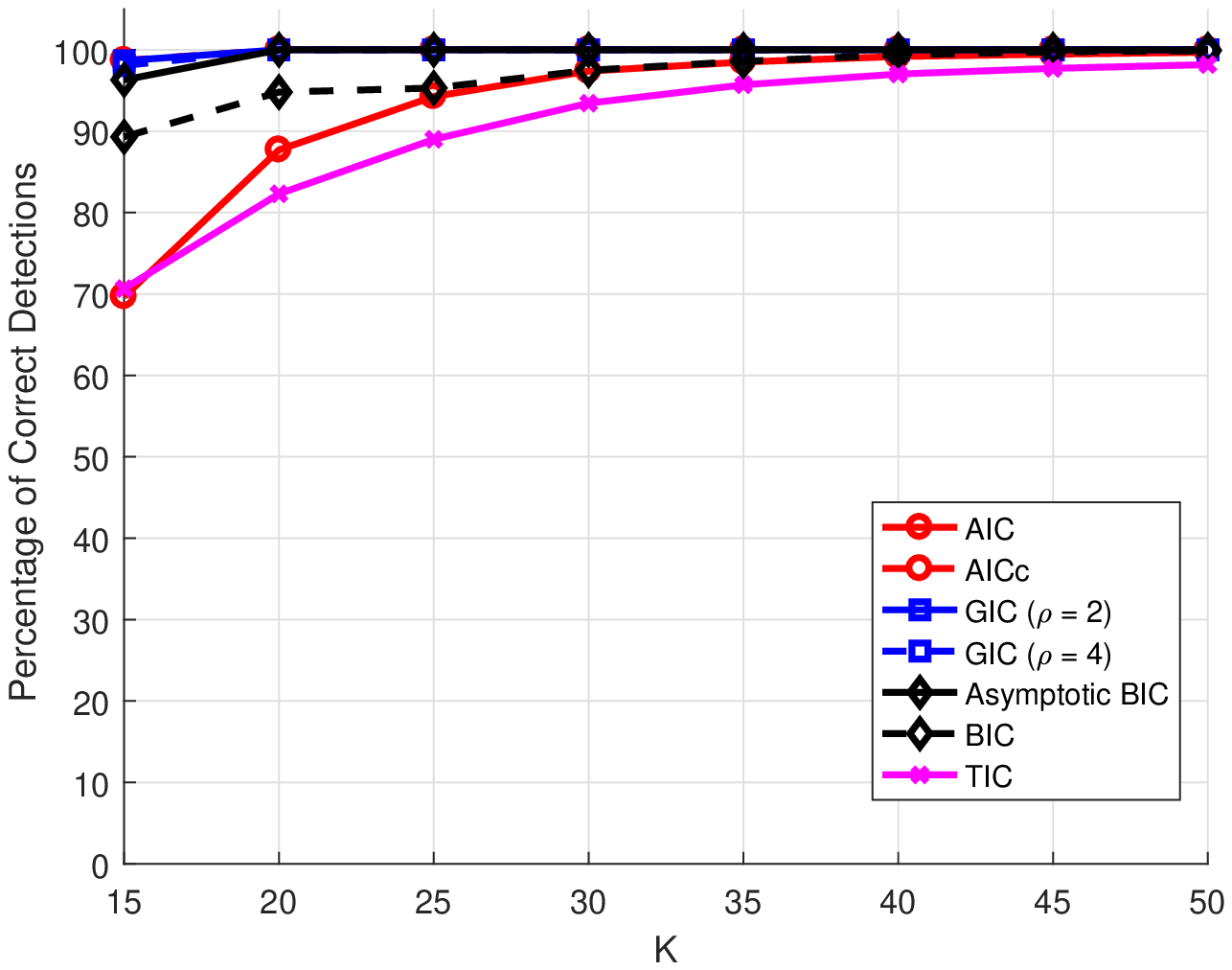}}
\subfigure[Hypothesis 4.]{\includegraphics[width=0.49\columnwidth]{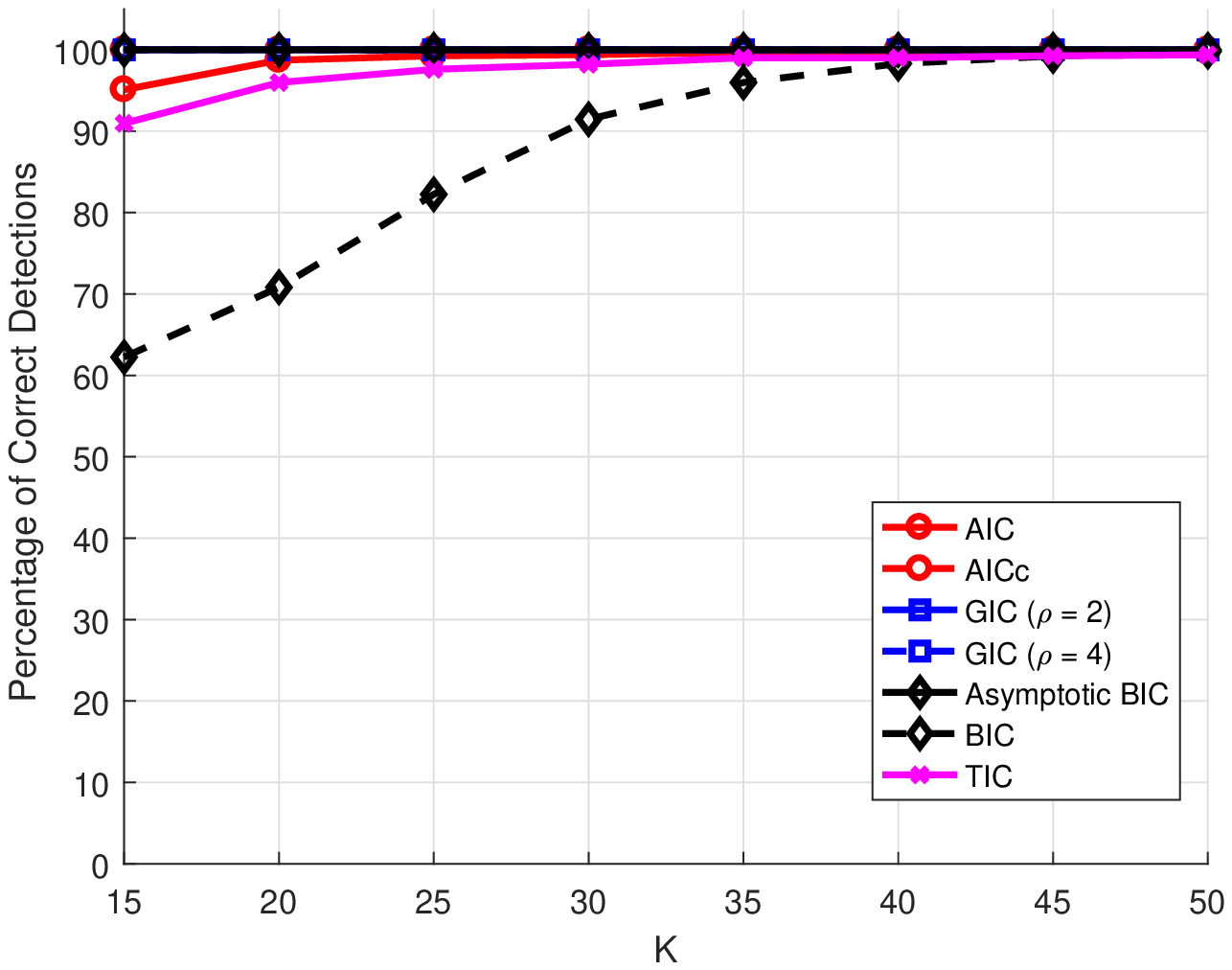}}
\caption{$P_{cc}$ versus $K$ for Study Case 1 and Approach A (primary and secondary data).}
\label{fig:Case1ApprA}
\end{figure}

\begin{figure}[H] \centering
\subfigure[Hypothesis 1.]{\includegraphics[width=0.49\columnwidth]{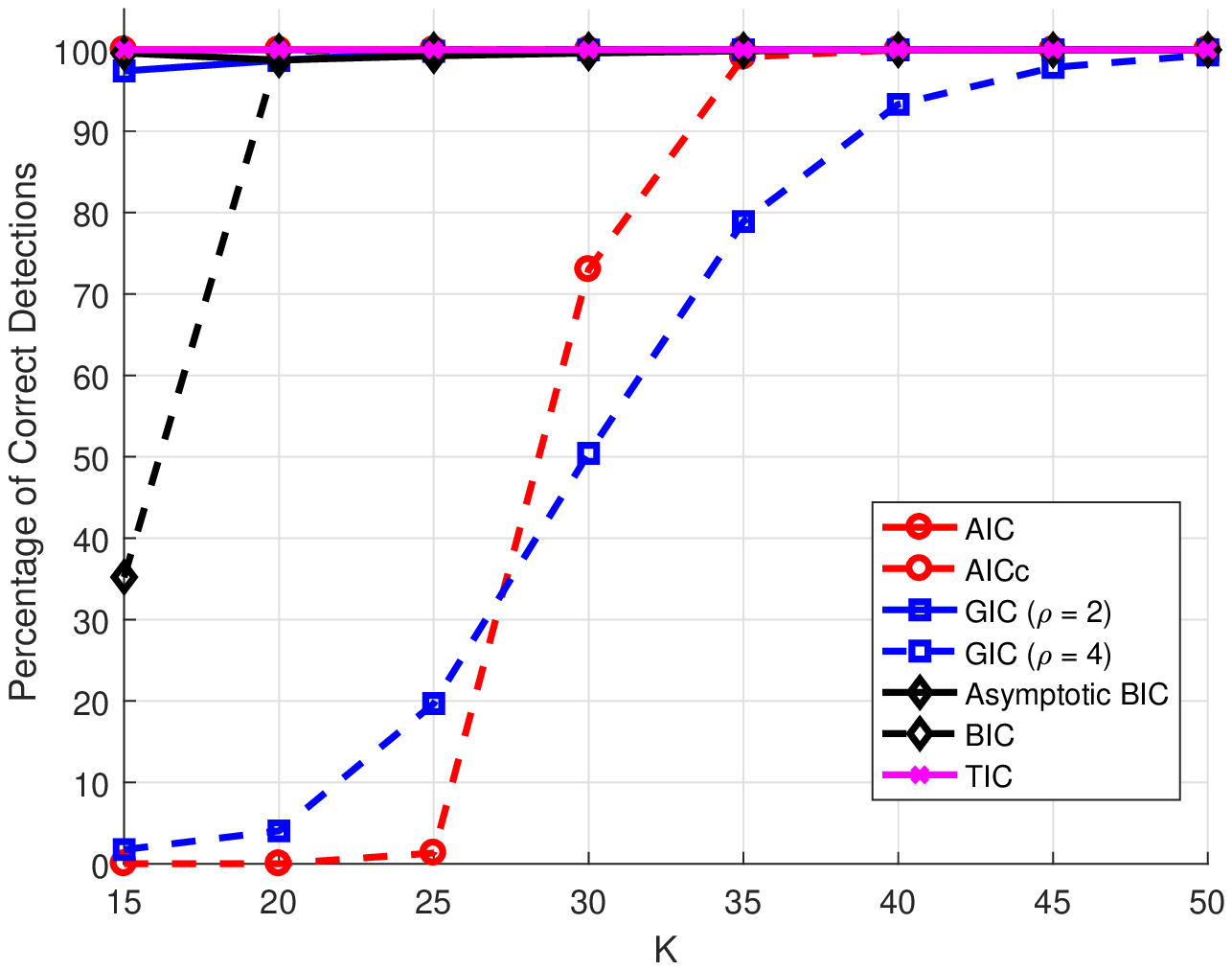}}
\subfigure[Hypothesis 2.]{\includegraphics[width=0.49\columnwidth]{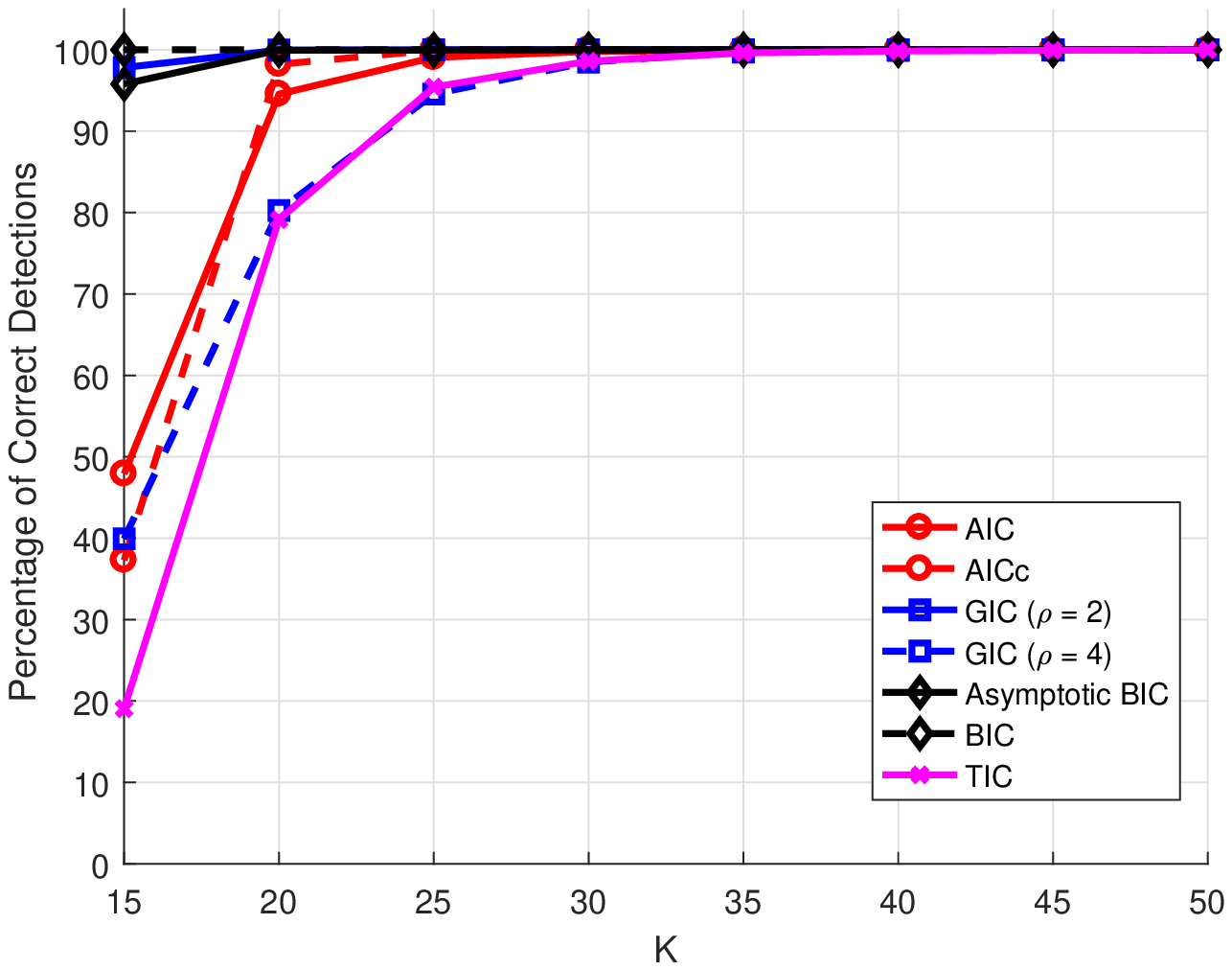}}\\
\subfigure[Hypothesis 3.]{\includegraphics[width=0.49\columnwidth]{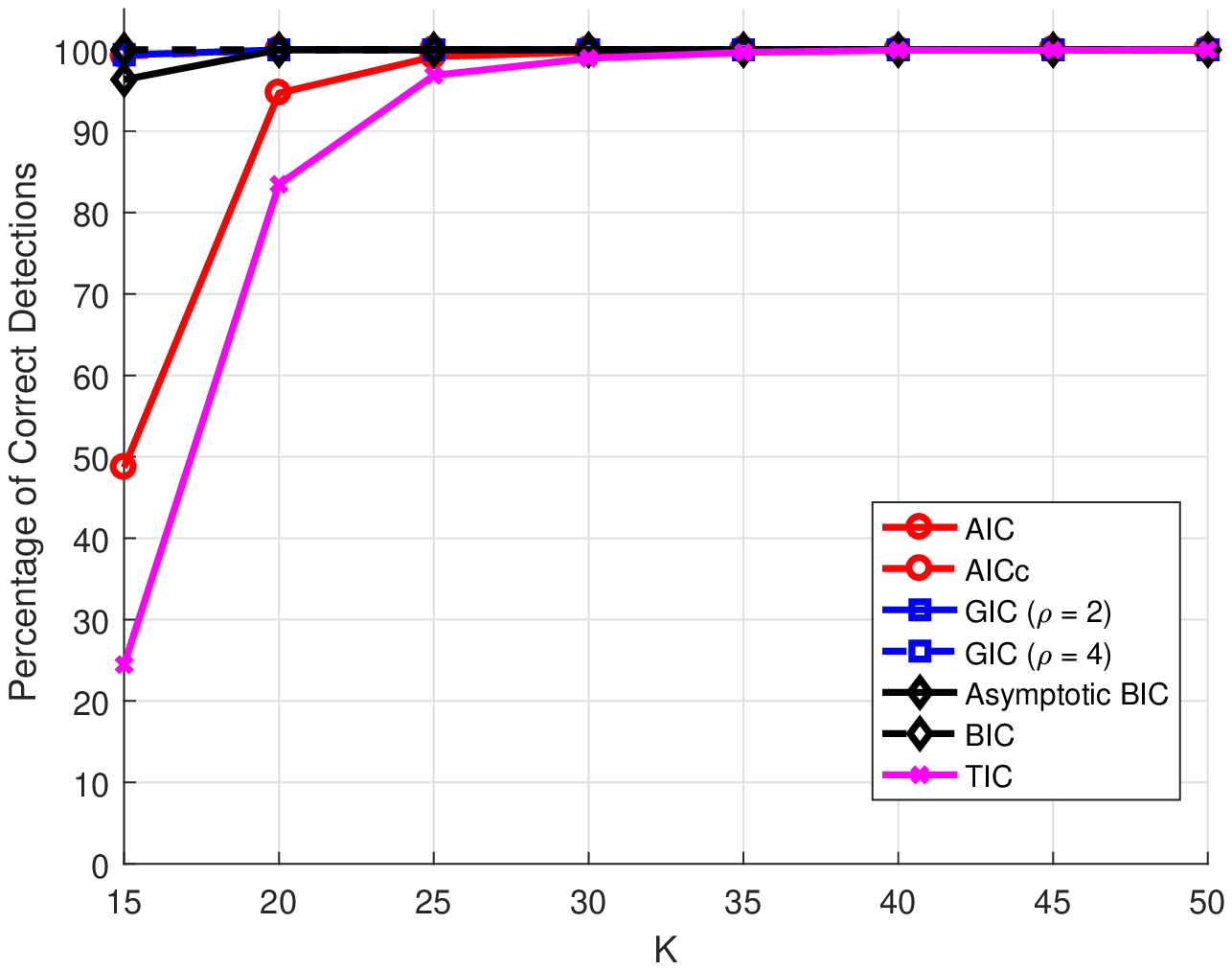}}
\subfigure[Hypothesis 4.]{\includegraphics[width=0.49\columnwidth]{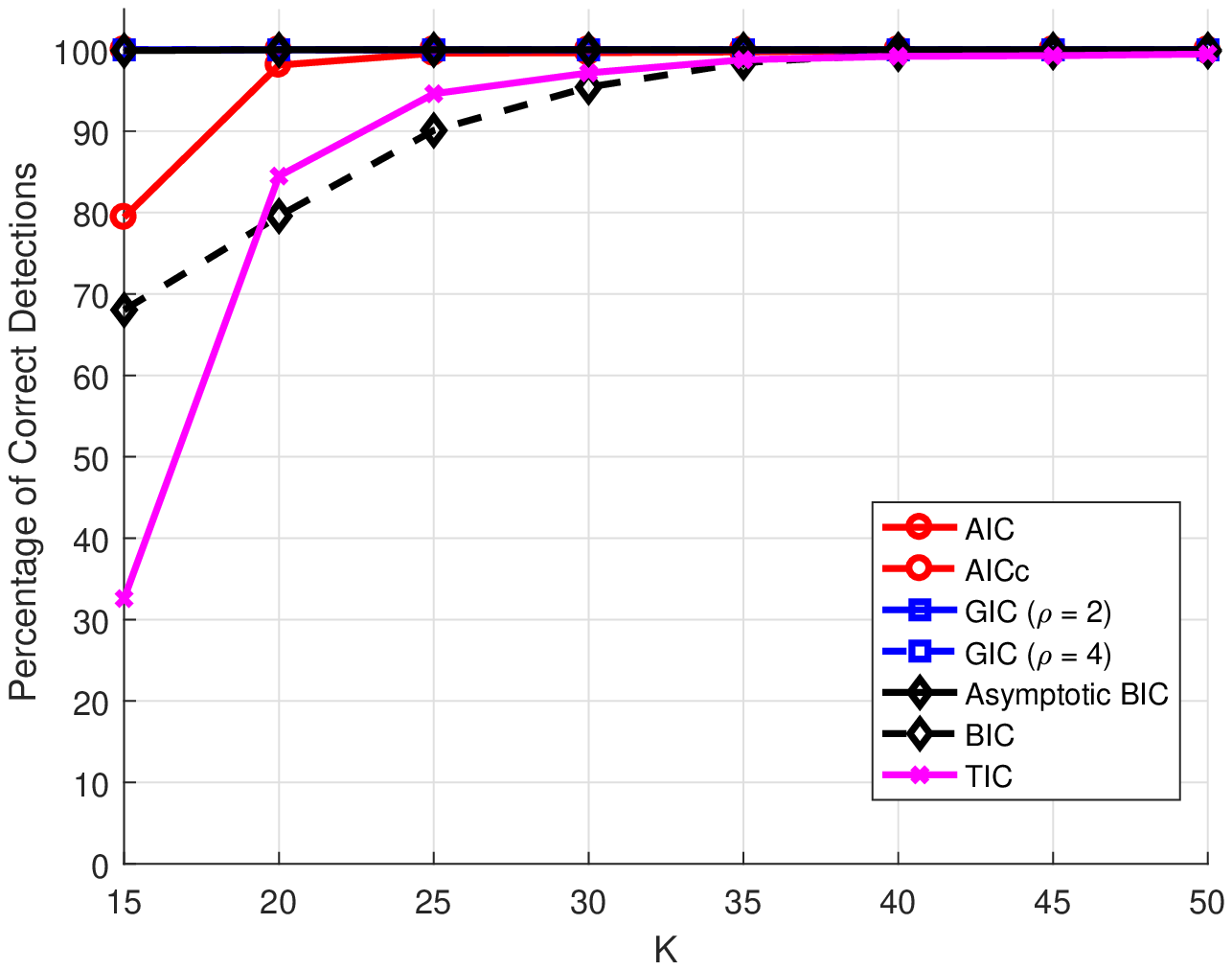}}
\caption{$P_{cc}$ versus $K$ for Study Case 1 and Approach B (secondary data only).}
\label{fig:Case1ApprB}
\end{figure}

\begin{figure}[H] \centering
\subfigure[Hypothesis 1.]{\includegraphics[width=0.49\columnwidth]{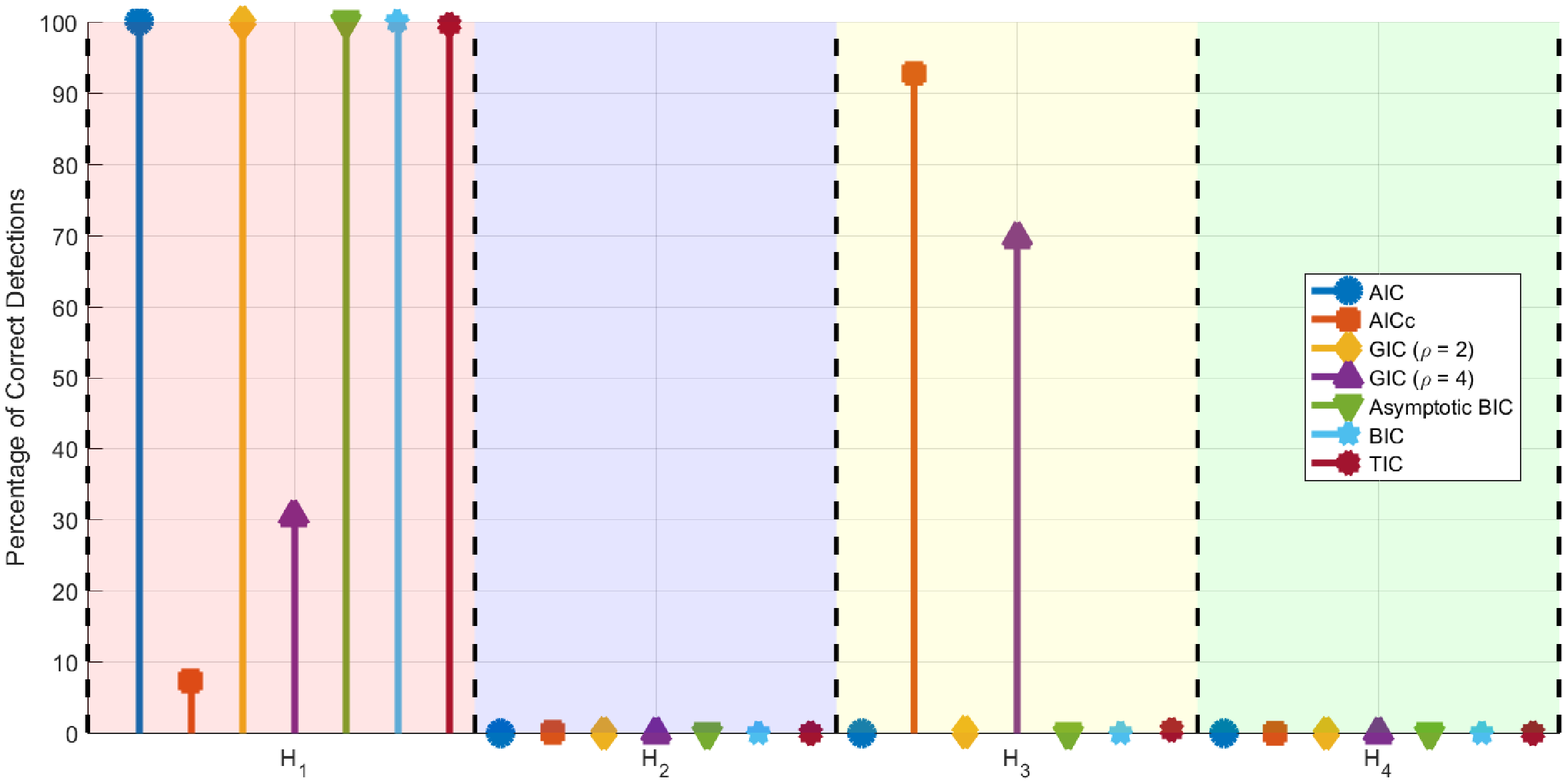}}
\subfigure[Hypothesis 2.]{\includegraphics[width=0.49\columnwidth]{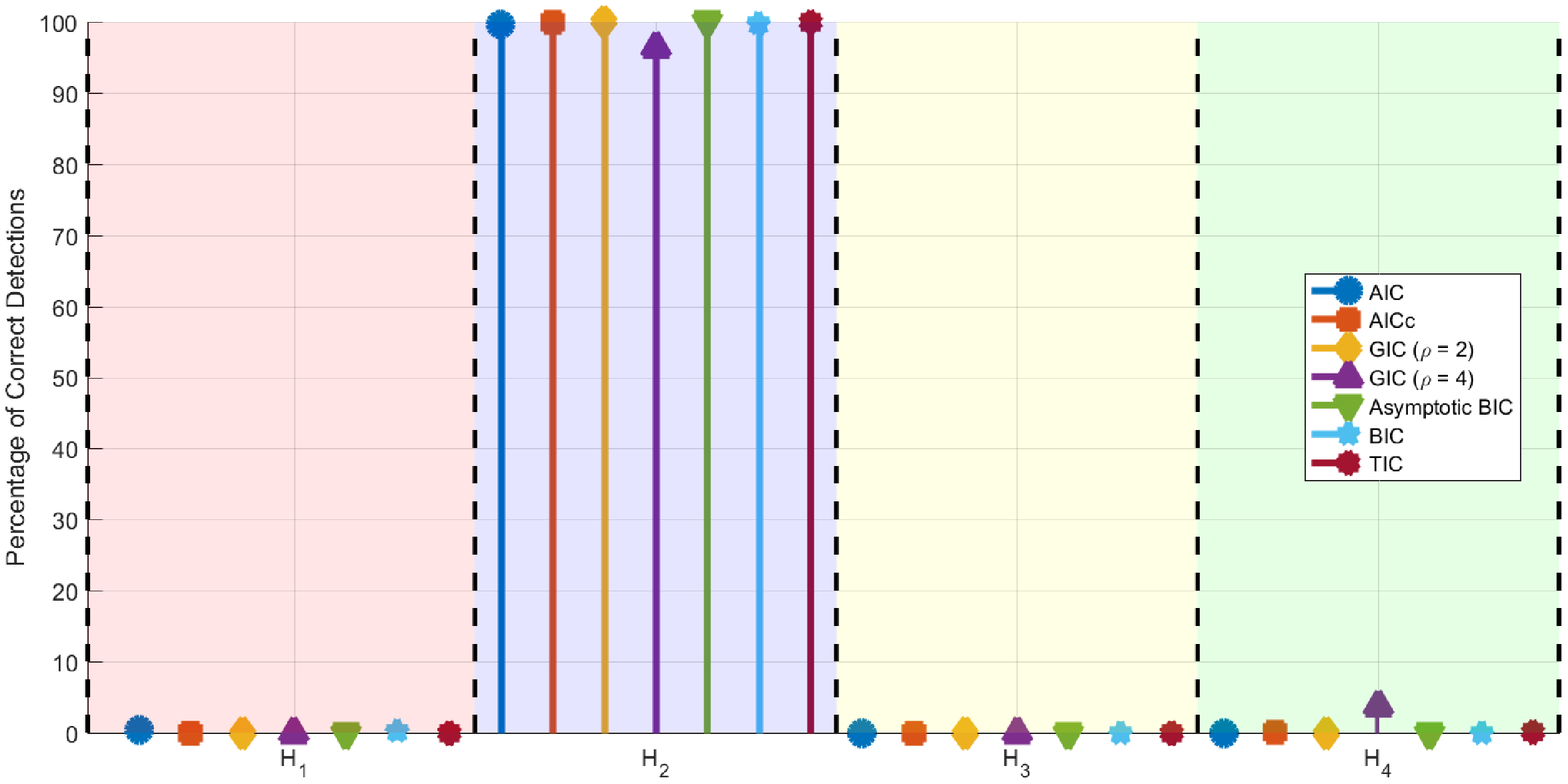}}\\
\subfigure[Hypothesis 3.]{\includegraphics[width=0.49\columnwidth]{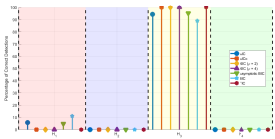}}
\subfigure[Hypothesis 4.]{\includegraphics[width=0.49\columnwidth]{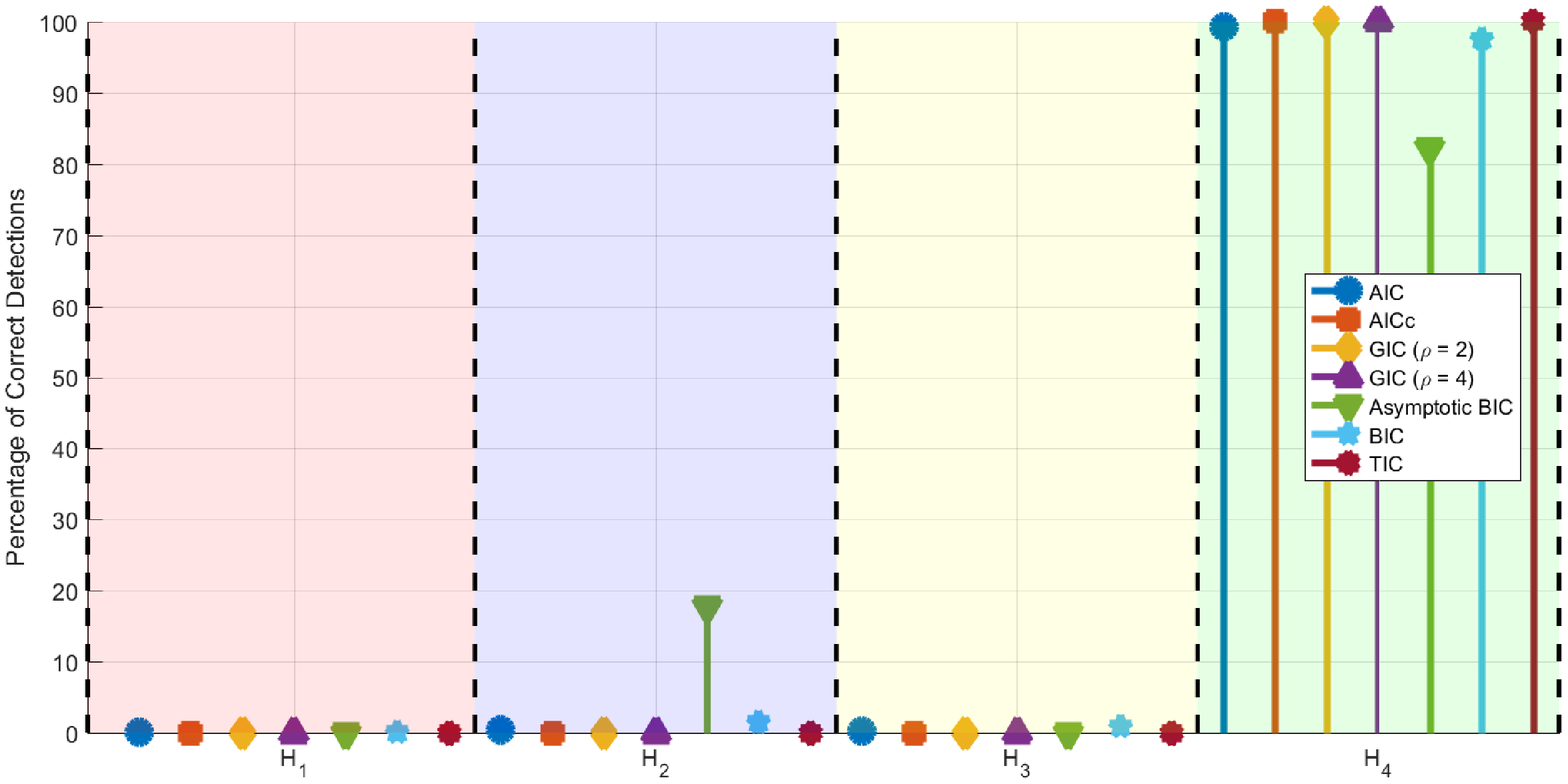}}
\caption{Percentage of classification for each hypothesis assuming Approach A and $K=25$.}
\label{fig:istCase1ApprA}
\end{figure}

Figures \ref{fig:Case1ApprA} and \ref{fig:Case1ApprB} refer to Case 1 and
contain the $P_{cc}$ curves for Approach A and B, respectively.
Inspection of the first figure highlights that AICc and GIC with $\rho=4$ exhibit poor performance under $H_1$ for $K<3N$ and under $H_2$ for $K<2N$.
This behavior is presumably due to the fact that in the current context AICc, as already stated, is heuristic, while the performance of GIC depends on the value of $\rho$.
Moreover, under $H_1$ and $H_4$, BIC requires $K>2N$ secondary data to achieve reasonable classification performances. Recall that BIC uses an estimate
of the FIM. The remaining classification schemes guarantee a $P_{cc}$ above $0.7$ over the considered range of values for $K$. The described trend remains the same
in Figure \ref{fig:Case1ApprB} except for a performance degradation for some architectures (such as AIC and TIC) when $K$ is low. 
The behavior of the considered rules can also be studied analyzing the classification percentages for each hypothesis. To this end,
in Figure \ref{fig:istCase1ApprA}, we plot the percentages of classification by means of histograms for Approach A and assuming $K=25$.
The inspection of the figure shows that under $H_1$ (or $H_2$), 
some MOS rules decide for $H_3$ (or $H_4$) and vice versa. In other words, the misclassification occurs between $H_1$ and $H_3$ or between
$H_2$ and $H_4$.
Finally, note that
including the CUT in the MOS classification rules (Approach A) leads to better performances than those obtained by means of Approach B.

\begin{figure}[H] \centering
\subfigure[Hypothesis 1.]{\includegraphics[width=0.49\columnwidth]{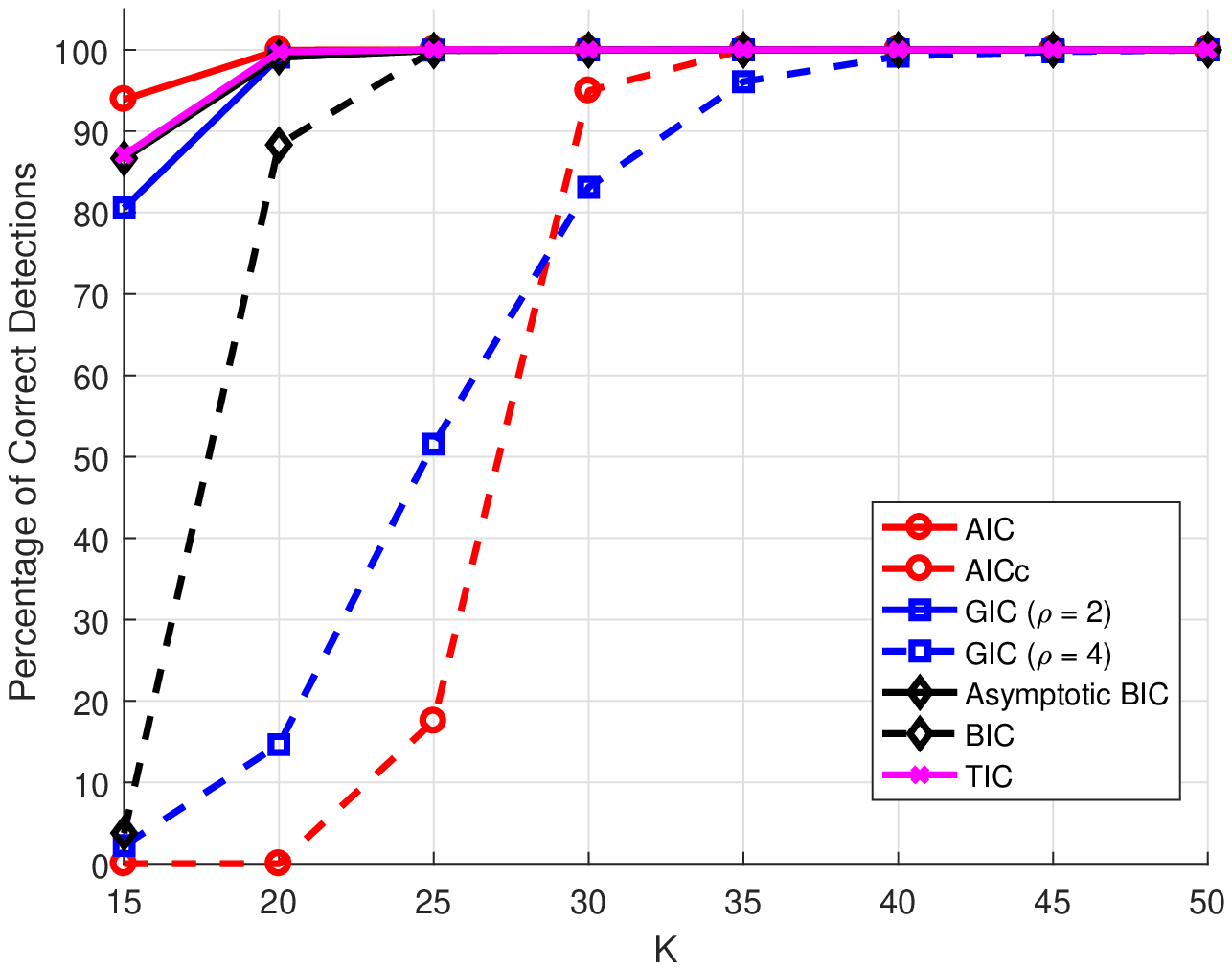}}
\subfigure[Hypothesis 2.]{\includegraphics[width=0.49\columnwidth]{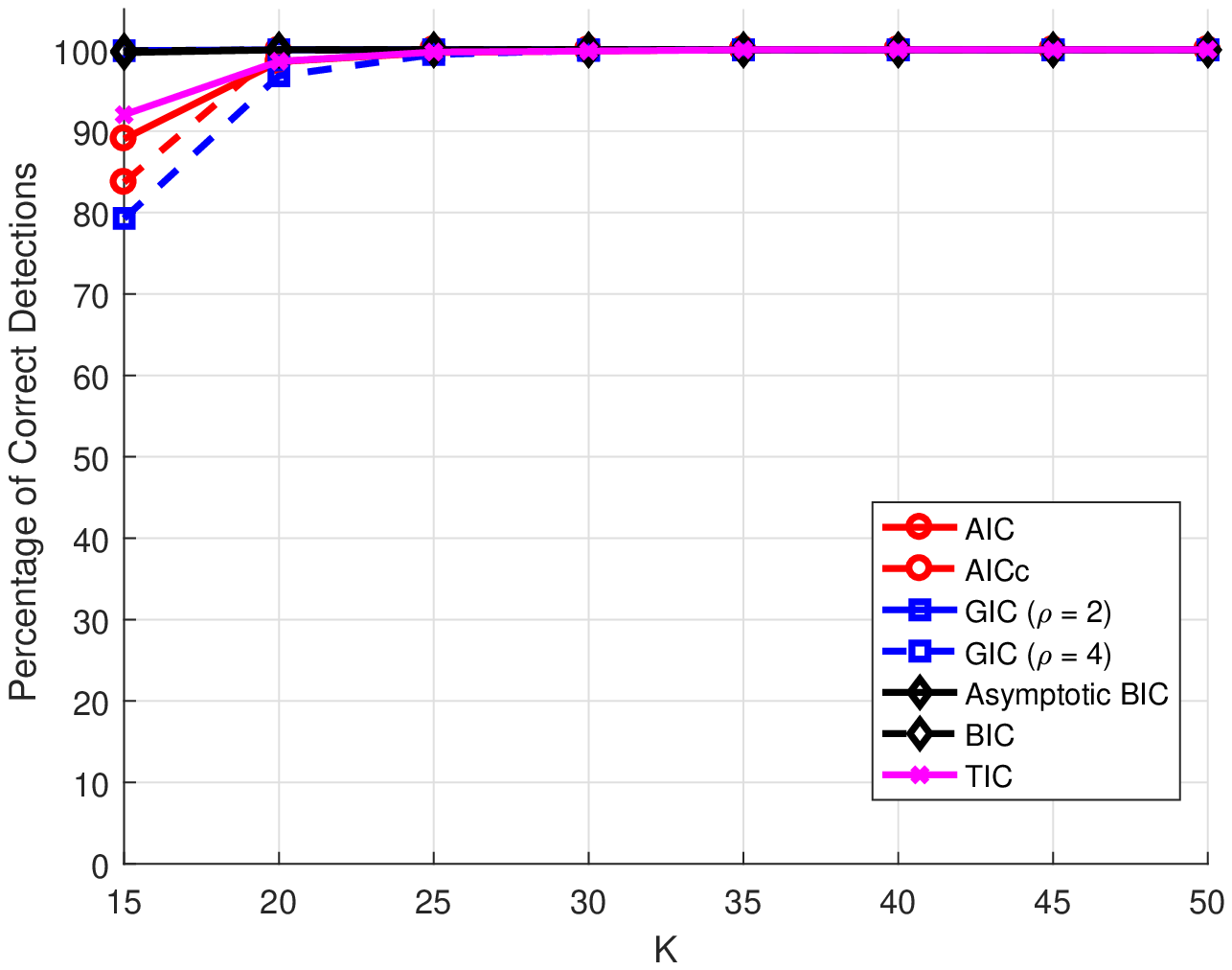}}\\
\subfigure[Hypothesis 3.]{\includegraphics[width=0.49\columnwidth]{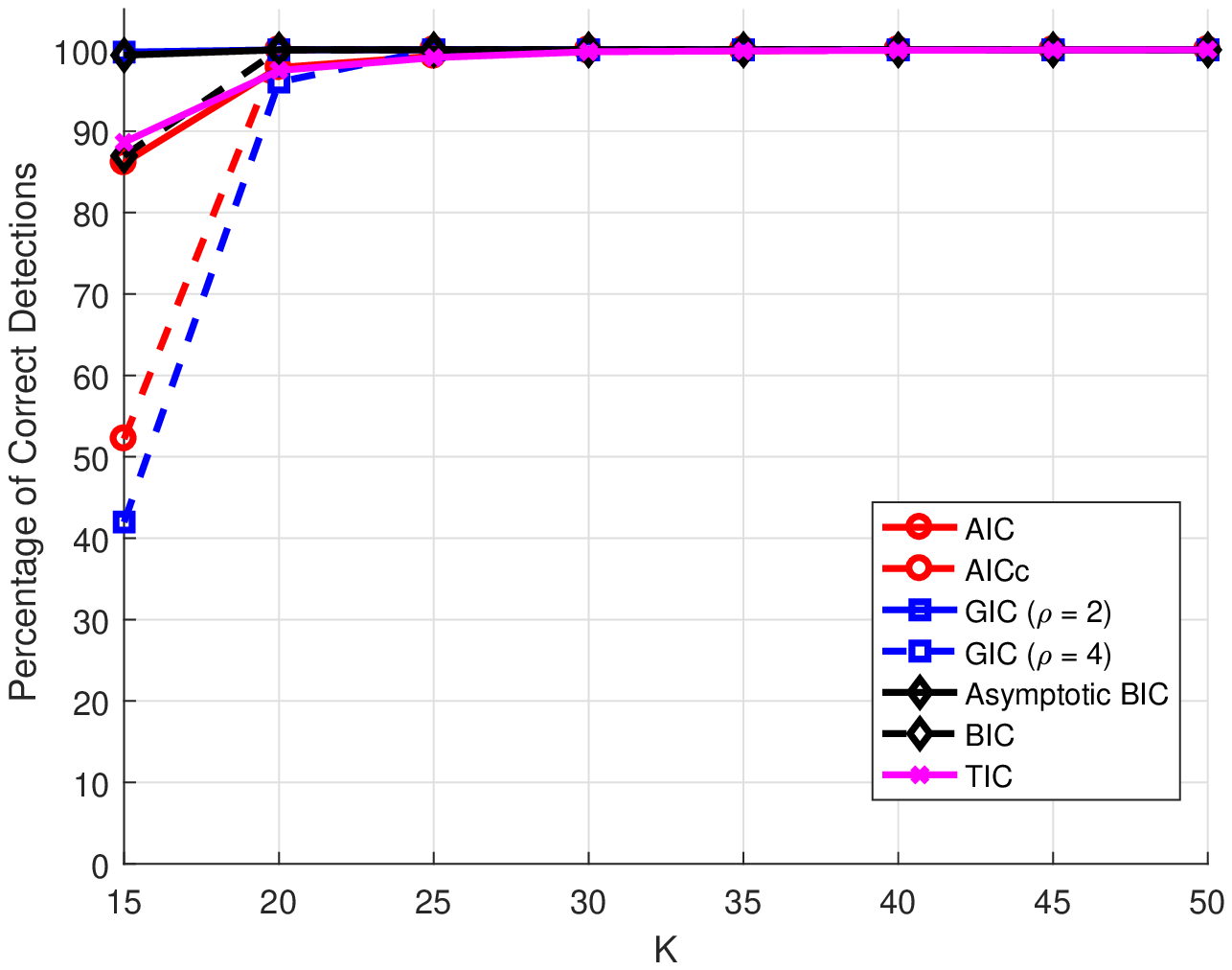}}
\subfigure[Hypothesis 4.]{\includegraphics[width=0.49\columnwidth]{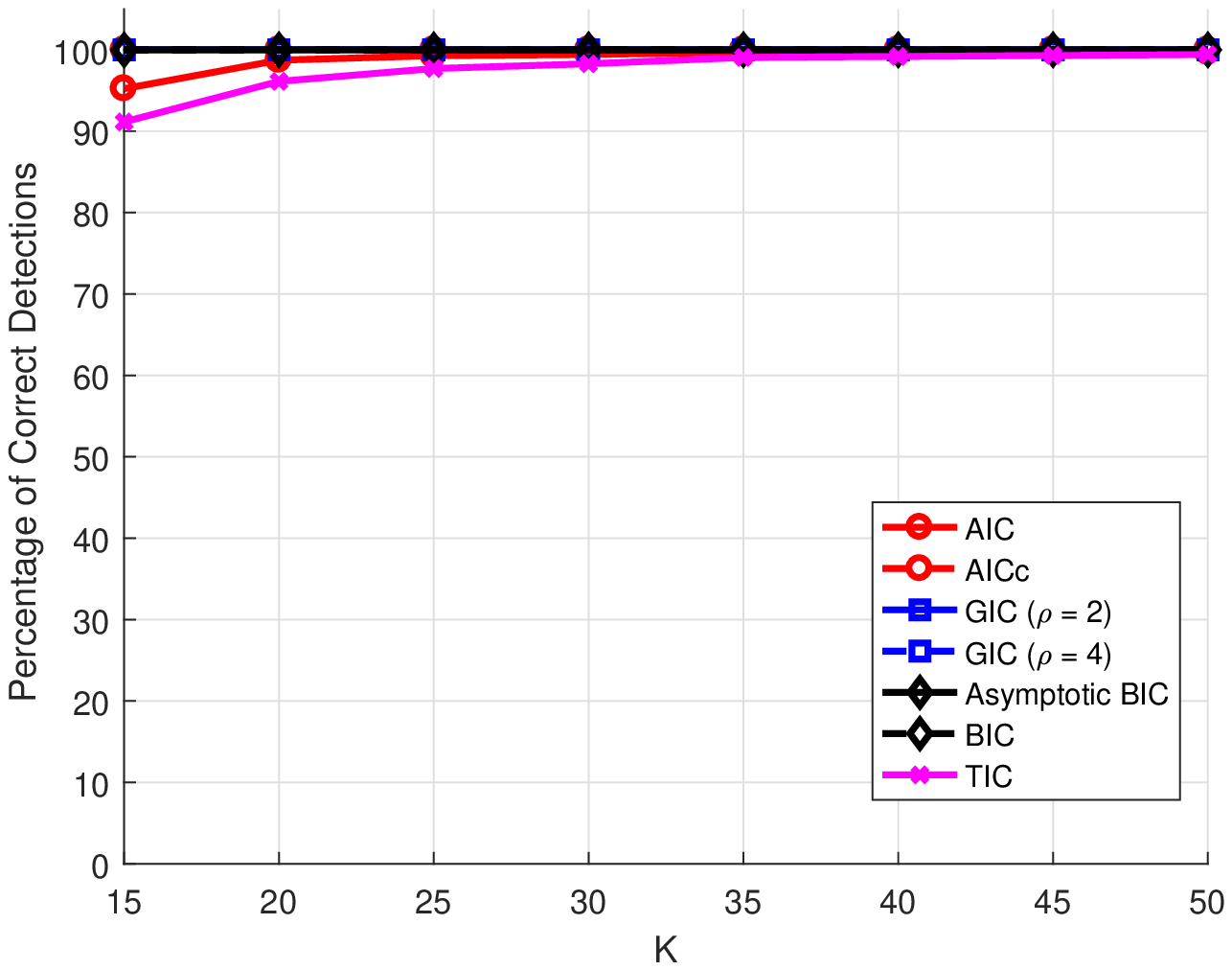}}
\caption{$P_{cc}$ versus $K$ for Study Case 2 and Approach A (primary and secondary data).}
\label{fig:Case2ApprA}
\end{figure}

\begin{figure}[H] \centering
\subfigure[Hypothesis 1.]{\includegraphics[width=0.49\columnwidth]{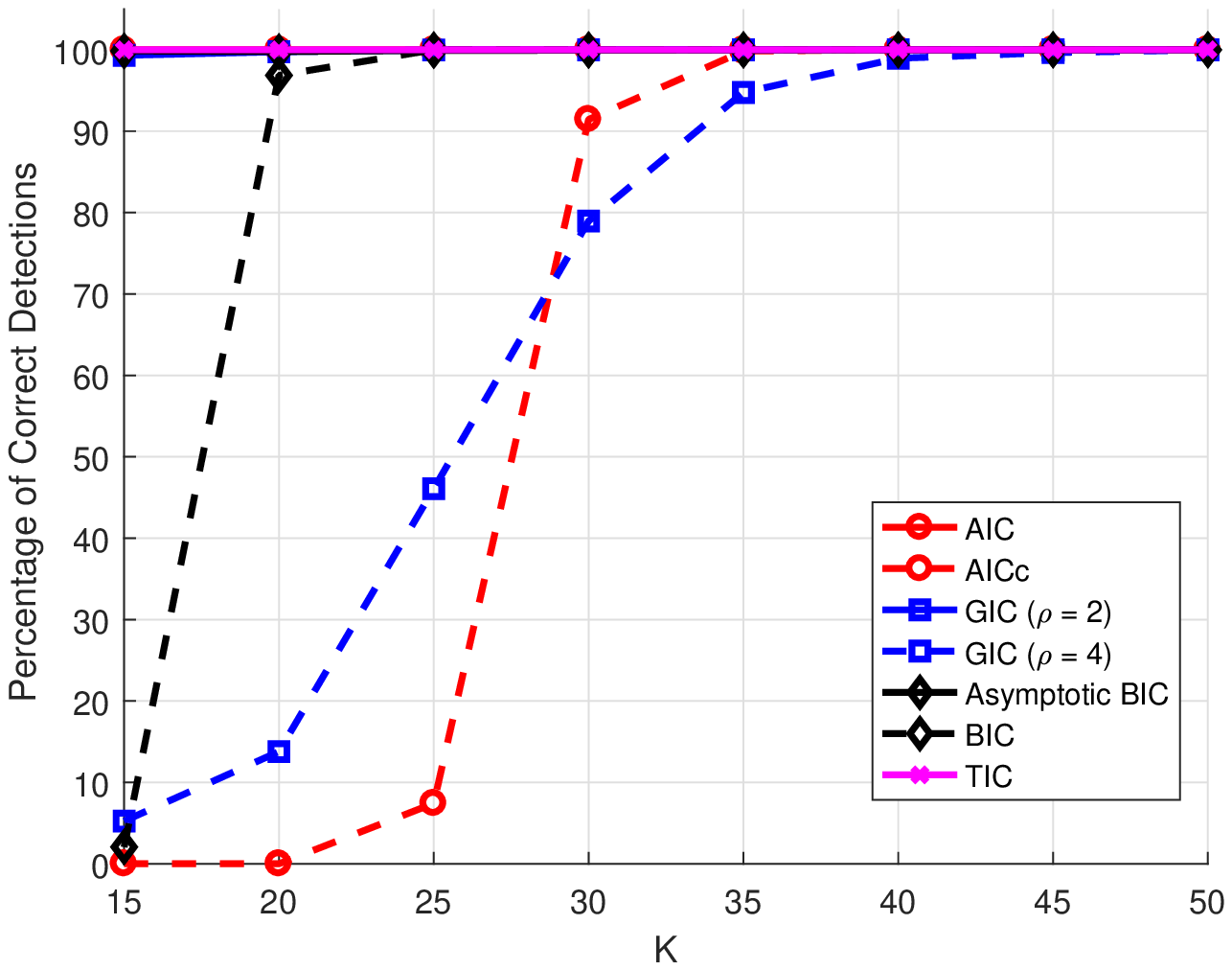}}
\subfigure[Hypothesis 2.]{\includegraphics[width=0.49\columnwidth]{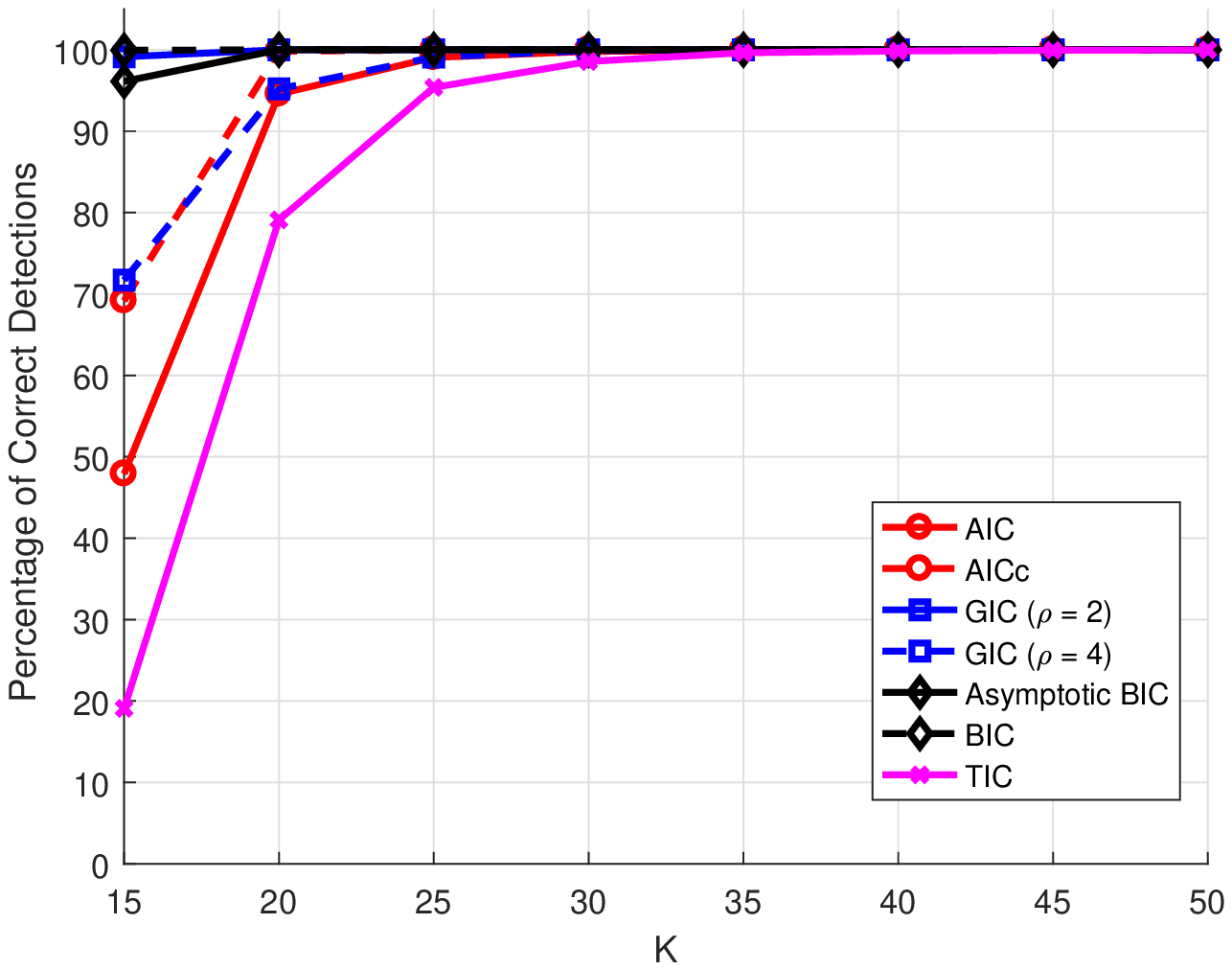}}\\
\subfigure[Hypothesis 3.]{\includegraphics[width=0.49\columnwidth]{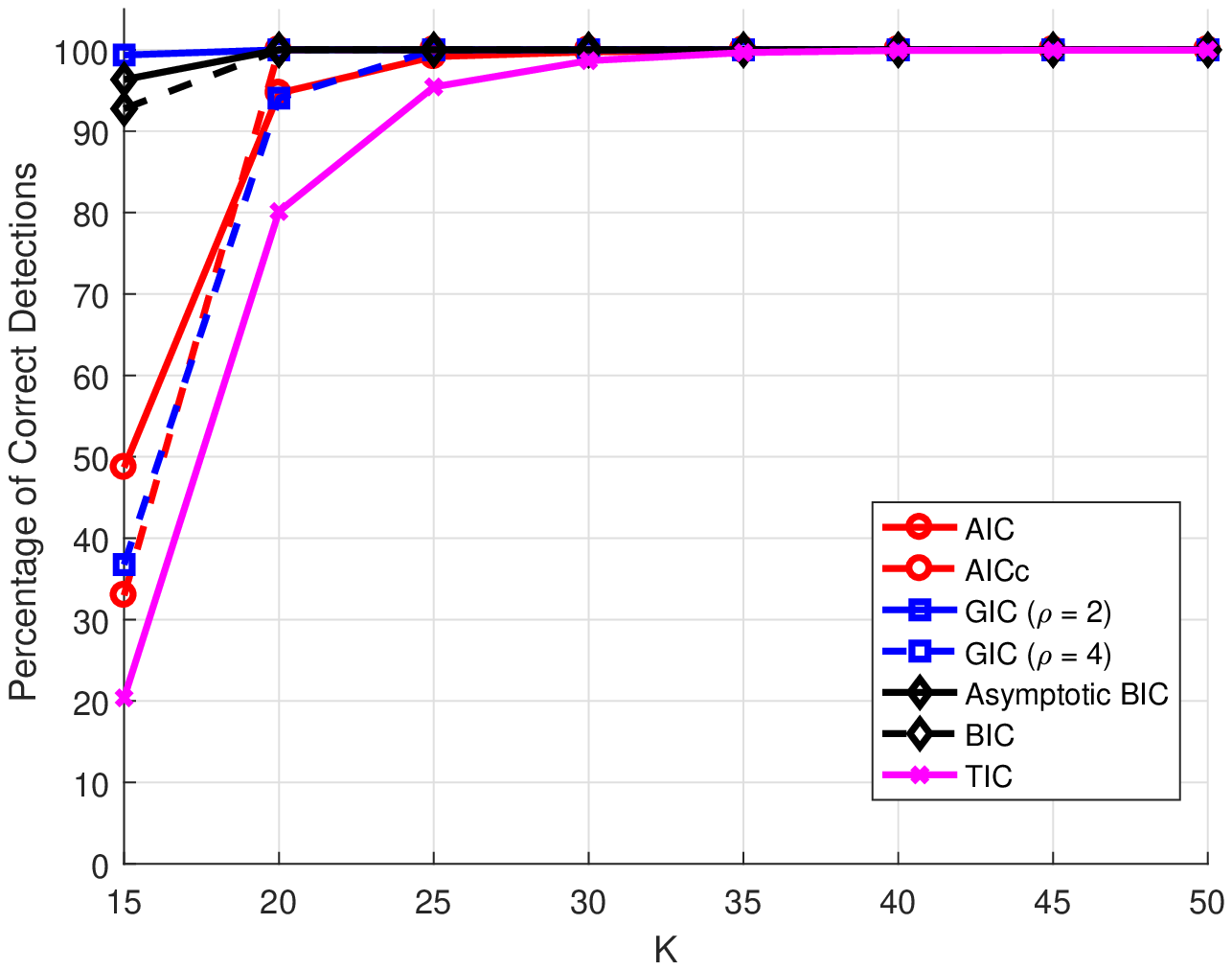}}
\subfigure[Hypothesis 4.]{\includegraphics[width=0.49\columnwidth]{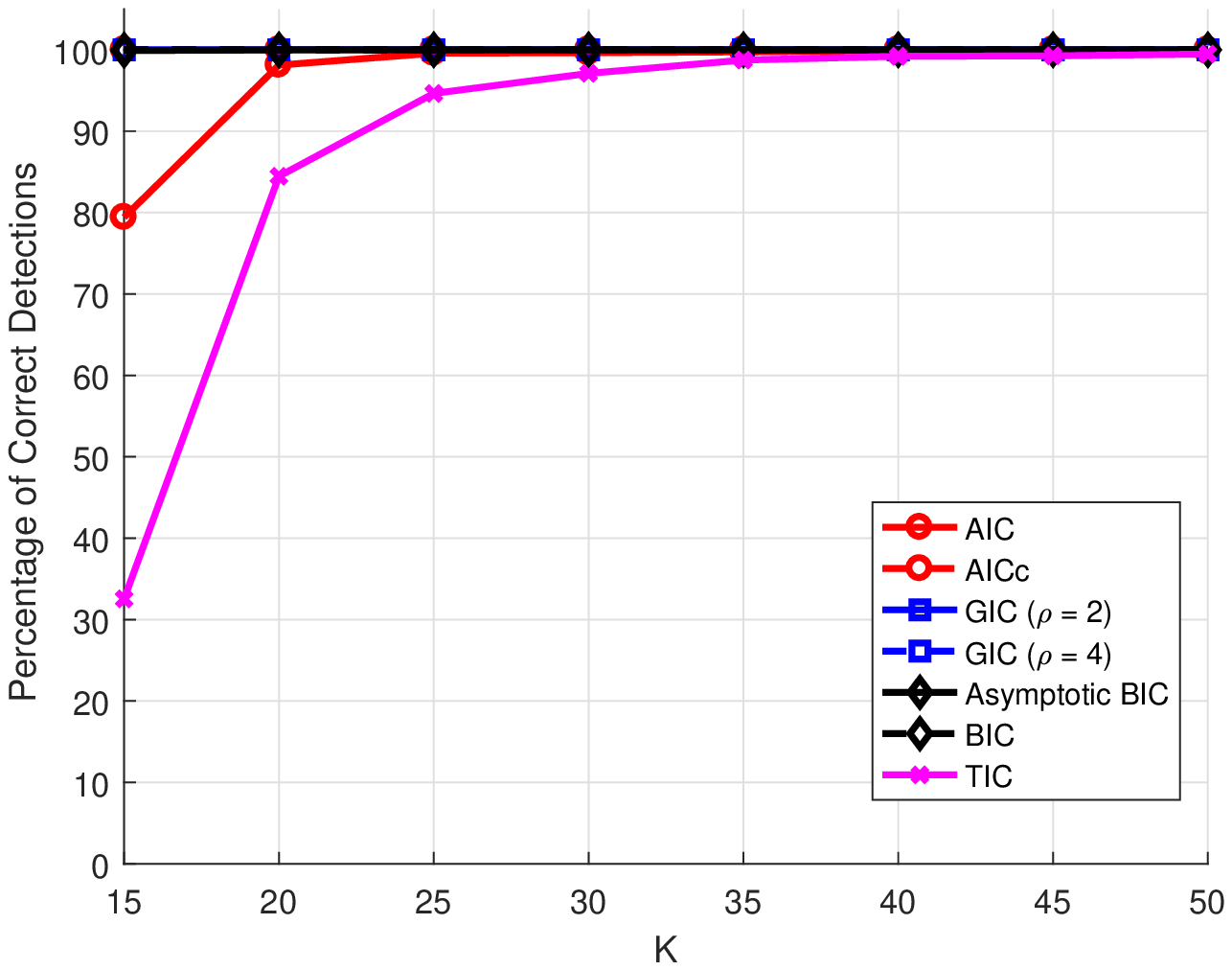}}
\caption{$P_{cc}$ versus $K$ for Study Case 2 and Approach B (secondary data only).}
\label{fig:Case2ApprB}
\end{figure}

In Figures \ref{fig:Case2ApprA} and \ref{fig:Case2ApprB}, the $P_{cc}$ curves for Case 2 are reported. The behavior of the classification rules is similar
to that observed in the previous figures with the difference that BIC suffers performance degradation for low values of $K$ under $H_1$ only.

From the inspection of all the above figures, it turns out that there does not exist a specific choice which provides the highest $P_{cc}$ under 
all the considered settings and parameters range. However,
the analysis underlines that the classification performances of some 
rules, in particular the AICc and GIC with $\rho=4$, are poor for low values of $K$ and this drawback
could be a reason to discard these architectures when $K\leq 2N$ and for the considered parameters setting.
In contrast to this, TIC and BIC classification schemes are
capable of guaranteeing $P_{cc}>0.8$ when $K\geq 2N$ in all 
the considered conditions. However, these rules become somewhat unstable when $K<2N$; this behavior may be due 
to the fact that the observed and sample FIM are less reliable when $K$ takes on relatively small values.
Finally, the Asymptotic BIC and GIC with $\rho=2$ provide the highest performance even for low values of $K$. 
The similarity in performance of these rules is due to the penalty terms whose values are close to each other 
for the considered parameters (i.e., $\log(K)\in [3, \ 3.8]$ for 
$K\in[20, \ 45]$). However, the hyperparameter $\rho$ of GIC is a degree of freedom that has to be suitably set 
(in fact, the GIC with $\rho=4$ has the worst performance), and there does not exist a general tuning criterion which allows us 
to choose the best value for $\rho$. 
On the other hand, the asymptotic BIC, which does not require any hyperparameter setting, stems as a reasonable 
operational choice at least for the considered scenarios.

\section{Conclusions}
\label{Sec:Conclusions}
This paper has considered the interference covariance structure classification which is of primary concern in some radar 
signal processing applications. Starting from a set of multivariate radar observations, the classification has been 
formulated as a multiple hypotheses test with some nested instances characterized by a different number of parameters. 
Several MOS rules, based on different theoretical criteria, have been devised to perform the covariance structure selection.
Besides, the possibilities of using primary and secondary data or only secondary vectors to implement the classification rules have been considered.
At the analysis stage their performance has been assessed in correspondence of two different operational scenarios highlighting the merits and the drawbacks connected with each approach. The classification curves, the complexity as well as the stability, has singled out the Asymptotic BIC based on 
secondary data only as the recommended selector for the considered scenarios.

Finally, two possible future research tracks deserve attention. First of all, 
we will study the effect of the proposed MOS techniques for ICM structure selection on the performance 
of target detection. Some preliminary results in this direction are encouraging: they show that using the 
proposed techniques leads to performances close to those of the oracle target detector that knows the actual structure of the ICM.
Then, the analysis on real radar data is essential to finally establish the effectiveness of the proposed approach.

\appendices

\section{Number of Parameters when $\bM$ is Centrohermitian or Centrosymmetric}
\label{App:PersymmetricPar}
Assume that $\bM\in\R^{N\times N}$ is centrosymmetric with $N$ even and let $m=N/2$; then $\bM$ can be partitioned as follows \cite{Reid}
\be
\bM=
\left[
\begin{array}{cc}
\bJ \bA \bJ & \bB^T
\\
\bB & \bA
\end{array}
\right],
\ee
where $\bA\in\R^{m \times m}$ is symmetric, $\bB\in\R^{m\times m}$ is persymmetric, 
and $\bJ$ is an $m$-dimensional permutation matrix. 
It is clear that
\begin{itemize}
\item the number of parameters defining $\bA$ is $m(m+1)/2$;
\item the number of parameters defining $\bB$ is $m(m+1)/2$.
\end{itemize}
Thus, $\bM$ can be represented by means of 
\be
m(m+1)=\frac{N}{2}\left(\frac{N}{2}+1\right)
\ee
parameters.

In the case where $N$ is still even and $\bM\in\C^{N\times N}$ is centrohermitian, $\bM$ has the following representation
\be
\bM=
\left[
\begin{array}{cc}
\bJ \bA^* \bJ & \bB^\dag
\\
\bB & \bA
\end{array}
\right],
\ee
where $\bA\in\C^{m \times m}$ is Hermitian and $\bB\in\C^{m\times m}$ persymmetric.
It follows that
\begin{itemize}
\item the number of parameters defining $\bA$ is $m^2$;
\item the number of parameters defining $\bB$ is $m(m+1)$.
\end{itemize}
The total number of parameters is 
\be
m(m+1) + m^2=\frac{N}{2}(N+1).
\ee

In order to complete the proof, assume that $N$ is odd and let $m=(N-1)/2$. 
Following the lead of \cite{Goldstein}, a centrosymmetric $\bM\in\R^{N \times N}$ can be partitioned as
\be
\bM=
\left[
\begin{array}{ccc}
\bJ \bA \bJ & \bc &\bB^T
\\
\bc^T & c & \bc^T \bJ
\\
\bB & \bJ \bc & \bA
\end{array}
\right],
\ee
where $\bA\in\R^{m \times m}$ is symmetric, $\bB\in\R^{m\times m}$ is persymmetric, $c\in\R$, and $\bc\in\R^{m\times 1}$.
It turns out that the total number of parameters is 
\be
m(m+1)+m+1=\left(\frac{N+1}{2}\right)^2.
\ee
Finally, assume that $\bM\in\C^{N\times N}$ is centrohermitian; then it can be partitioned as \cite{Goldstein}
\be
\bM=
\left[
\begin{array}{ccc}
\bJ \bA^* \bJ & \bc &\bB^\dag
\\
\bc^\dag & c & \bc^\dag \bJ
\\
\bB & \bJ \bc & \bA
\end{array}
\right],
\ee
where $\bA\in\C^{m \times m}$ is Hermitian, $\bB\in\C^{m\times m}$ is persymmetric, $c\in\R$, and $\bc\in\C^{m\times 1}$.
As a consequence, the number of parameters characterizing $\bM$ is 
\be
2m^2+3m+1=N(N+1)/2.
\ee

\section{Gradient of the Log-Likelihood Functions}
\label{App:Grad_z_zk}
As a preliminary remark, observe that the ICM is always either Hermitian or symmetric. Let us first 
focus on $s(\bp_i,H_i; \bz)$ and evaluate the first derivative of this function
with respect to the $l$th component of $\bp_{i}$. It follows that two cases are possible: $\bp_i(l)$ is a component of $\btheta_i$ or
$\bp_i(l)$ is a component of $\balpha$.

As for the first case, it is possible to show that
\begin{align}
\frac{\partial s(\bp_{i}, H_i ; \bz)}{\partial \btheta_i(l)}
&=- \frac{\partial}{\partial \btheta_{i}(l)} \left\{
\log\det(\bM_i)
\right\}-
\frac{\partial}{\partial \btheta_{i}(l)} \left\{
\tr\left[\bX_i\bS_{\alpha}\right]
\right\}\nonumber
\\
&=-\left\{\vect\left[(\bX_i)^T\right]\right\}^T\frac{\partial }{\partial\btheta_{i}(l)}\{\vect [\bM_i]\} \nonumber
\\
& +\tr\left\{ \bX_i \bS_{\alpha} \bX_i  \frac{\partial \bM_i}{\partial\btheta_{i}(l)}  \right\},
\label{eqn:firstDerShom}
\end{align}
where the last equality comes from equations $A.390$ and $A.391$ of \cite{VanTrees4}. The above equation can be further simplified observing that
\be
\vect[\bM_i]=\bC_i\btheta_{i},
\label{eqn:vecMC}
\ee
where $\bC_i \in \C^{N^2 \times m_i}$ is a transformation matrix that depends on the specific structure of 
$\bM_i$ and on how $\btheta_{i}$ is defined. For instance, if $\bM_i$ is Hermitian unstructured with $N=3$ and
\begin{equation}
\btheta_{i} = 
\begin{bmatrix}
\bM(1,1)\\
\Re \{ \bM(2,1) \}\\
\Im \{ \bM(2,1) \}\\
\Re \{ \bM(3,1) \}\\
\Im \{ \bM(3,1) \}\\
\bM(2,2)\\
\Re \{ \bM(3,2) \}\\
\Im \{ \bM(3,2) \}\\
\bM(3,3)\\
\end{bmatrix},
\end{equation}
then 
\begin{equation}
\bC_i = 
\begin{bmatrix}
1 & 0 & 0 & 0 & 0 & 0 & 0 & 0 & 0\\
0 & 1 & j & 0 & 0 & 0 & 0 & 0 & 0\\
0 & 0 & 0 & 1 & j & 0 & 0 & 0 & 0\\
0 & 1 & -j & 0 & 0 & 0 & 0 & 0 & 0\\
0 & 0 & 0 & 0 & 0 & 1 & 0 & 0 & 0\\
0 & 0 & 0 & 0 & 0 & 0 & 1 & j & 0\\
0 & 0 & 0 & 1 & -j & 0 & 0 & 0 & 0\\
0 & 0 & 0 & 0 & 0 & 0 & 1 & -j & 0\\
0 & 0 & 0 & 0 & 0 & 0 & 0 & 0 & 1\\
\end{bmatrix}.
\end{equation}
It follows that
\be
\frac{\partial }{\partial \btheta_{i}(l)}[\bC_i\btheta_{i}]=\bC_i \boe_{l,i},
\label{eqn:vecC}
\ee
where $\boe_{l,i}$ is the $l$th elementary vector of size $m_i$. 
Moreover, let $\bA$, $\bB$, $\bC$, and $\bD$ be generic matrices whose sizes are such that the product
$\bA\bB\bC\bD$ makes sense and yields a square matrix; then the following equality holds \cite{Magnus}
\be
\tr(\bA\bB\bC\bD)=[\vect(\bD^T)]^T(\bC^T \otimes \bA) \vect(\bB).
\label{eqn:traceProduct}
\ee
Thus, the second term of \eqref{eqn:firstDerShom} can be recast as
\begin{multline}
\tr\left\{ \bX_i \bS_{\alpha} \bX_i  \frac{\partial \bM_i}{\partial\btheta_{i}(l)}  \right\}
\\
=\left\{\vect\left[ \frac{\partial \bM^T(\btheta_{i})}{\partial \btheta_{i}(l)} \right]\right\}^T
[(\bX_i)^T\otimes \bX_i] \vect[\bS_{\alpha}]
\end{multline}
Gathering the above results and accounting for $\bM_i$ being symmetric or Hermitian, \eqref{eqn:firstDerShom} becomes
\be
\mbox{\eqref{eqn:firstDerShom}}=
\begin{cases}
-\left\{\vect\left[\bX_i^*\right]\right\}^T\bC_i\boe_{l,i}+\boe_{l,i}^T \bC_i^\dag \times \\ 
[\bX_i^*\otimes \bX_i] \vect[\bS_{\alpha}], \mbox{ if } \bM_i \mbox{ is Hermitian},
\\
-\left\{\vect\left[(\bX_i)\right]\right\}^T\bC_i\boe_{l,i}+\boe_{l,i}^T \bC_i^T \times\\
[\bX_i\otimes \bX_i] \vect[\bS_{\alpha}], \mbox{ if } \bM_i \mbox{ is symmetric},
\end{cases}
\label{eqn:derivativeFinalHom}
\ee
where the following equality has been used
\begin{equation} \label{eqn:vedDerC}
\begin{split}
& \vect\left[ \frac{\partial \bM^T(\btheta_{i})}{\partial \btheta_{i}(l)} \right]=\\
&
\begin{cases}
\ds\vect\left[ \frac{\partial \bM^*(\btheta_{i})}{\partial \btheta_{i}(l)} \right]=
\frac{\partial }{\partial \btheta_{i}(l)}\{[\vect(\bM_i)]^*\}=\\
\qquad \bC^*_i\boe_{l,i}, \mbox{ if } \bM_i \mbox{ is Hermitian},
\\
\vspace{-5mm}
\\
\ds\vect\left[ \frac{\partial \bM_i}{\partial \btheta_{i}(l)} \right]=
\frac{\partial }{\partial \btheta_{i}(l)}[\bC_i\btheta_{i}] = \\
\qquad \bC_i\boe_{l,i}, \mbox{ if } \bM_i \mbox{ is symmetric}.
\end{cases}
\end{split}
\end{equation}
Hence, exploiting \eqref{eqn:derivativeFinalHom}, it is not difficult to obtain \eqref{eqn:gradient_z}. 
Following the same line of reasoning and replacing $\bS_\alpha$ with $\bS_k$, it is possible to prove \eqref{eqn:gradient_zk}.

As a final step, we evaluate the gradient of $s(\bp_i,H_i;\bz)$ with respect to $\balpha$. To this end, observe that
\begin{equation}
\begin{split}
& \frac{\partial s(\bp_i,H_i;\bz)}{\partial \balpha}= \frac{\partial}{\partial \balpha} \left\{-\tr[\bX_i\bS_\alpha]  \right\} =\\
& \frac{\partial}{\partial \balpha} \left\{  
\alpha \bz^\dag \bX_i \bv + \alpha^* \bv^\dag \bX_i \bz -
\alpha \alpha^* \bv^\dag \bX_i \bv \right\} .\\
\end{split}
\end{equation}
Using the above equation, the gradient with respect to $\balpha$ can be expressed as in \eqref{eqn:gradient_zAalpha}.

\section{Hessian of the Log-Likelihood Function}
\label{App:Hess}

In this appendix we derive the Hessian of the log-likelihood function $s(\bp_i,H_i;\bz,\bZ)$. 
To this end, consider $\bH_{\theta\theta,i}$, whose $(l,m)$-entry can be written as
\begin{equation} \label{eqn:Hess00}
\begin{split}
&\bH_{\theta\theta,i}(l,m)=\frac{\partial^2 s(\bp_{i}, H_i ;\bz, \bZ)}{\partial \btheta_i(l) \partial \btheta_i(m)} = \\
& -(K+1) \frac{\partial^2 }{\partial \btheta_{i}(l)\partial \btheta_{i}(m)} [\log\det(\bM_i)]\\
& -\frac{\partial^2 }{\partial \btheta_{i}(l)\partial \btheta_{i}(m)} [\tr(\bX_i(\bS_{\alpha}+\bS))]
\\
&=(K+1)\tr\left\{ \bX_i \frac{\partial \bM_i}{\partial \btheta_{i}(l)} 
\bX_i \frac{\partial \bM_i}{\partial \btheta_{i}(m)} \right\}\\ 
& - (K+1) \tr\left\{ \bX_i \frac{\partial^2 \bM_i}{\partial \btheta_{i}(l)\partial \btheta_{i}(m)} \right\}\\
& + \frac{\partial }{\partial \btheta_{i}(m)} \left\{
\tr\left\{
\bX_i (\bS_{\alpha}+\bS) \bX_i \frac{\partial \bM_i}{\partial \btheta_{i}(l)}
\right\}
\right\},
\end{split}
\end{equation}
where the last equality comes from the application of $(A.391)$ and $(A.393)$ in \cite{VanTrees4}. Now, let us focus on the last term of \eqref{eqn:Hess00} and exploit $(A.391)$ of \cite{VanTrees4} to obtain
\begin{align}
&\frac{\partial }{\partial \btheta_{i}(m)} \left\{
\tr\left\{
\bX_i (\bS_{\alpha}+\bS) \bX_i \frac{\partial \bM_i}{\partial \btheta_{i}(l)}
\right\}
\right\}\nonumber
\\
&=\tr\left\{
\frac{\partial }{\partial \btheta_{i}(m)} \left\{
\bX_i (\bS_{\alpha}+\bS) \bX_i \frac{\partial \bM_i}{\partial \btheta_{i}(l)}
\right\}
\right\}\nonumber
\\
&=-\tr\left\{
(\bS_{\alpha}+\bS) \bX_i \frac{\partial \bM_i}{\partial \btheta_{i}(l)} \bX_i 
\frac{\partial \bM_i}{\partial \btheta_{i}(m)} \bX_i
\right\}\nonumber
\\
& + \tr \Big\{ \bX_i (\bS_{\alpha}+\bS) \times \nonumber
\\
& \left[
\frac{\partial \bX_i}{\partial \btheta_{i}(m)}\frac{\partial \bM_i}{\partial \btheta_{i}(l)}
+\bX_i \frac{\partial^2 \bM_i}{\partial \btheta_{i}(l)\partial \btheta_{i}(m)}
\right] \Big\} \nonumber
\\
&=-\tr\left\{
(\bS_{\alpha}+\bS)\bX_i \frac{\partial \bM_i}{\partial \btheta_{i}(l)} \bX_i 
\frac{\partial \bM_i}{\partial \btheta_{i}(m)} \bX_i
\right\}\nonumber
\\
&
-\tr\left\{
(\bS_{\alpha}+\bS) \bX_i \frac{\partial \bM_i}{\partial \btheta_{i}(m)} \bX_i 
\frac{\partial \bM_i}{\partial \btheta_{i}(l)} \bX_i
\right\}\nonumber
\\
&+\tr\left\{
\bX_i (\bS_{\alpha}+\bS) \bX_i \frac{\partial^2 \bM_i}{\partial \btheta_{i}(l)\partial \btheta_{i}(m)}
\right\}.
\end{align}
The terms involving the second-order derivative of $\bM_i$ can be discarded because
\begin{equation}
\begin{split}
\vect\left[
\frac{\partial^2 \bM_i}{\partial \btheta_{i}(l)\partial \btheta_{i}(m)}
\right]
& =\frac{\partial^2 \vect\left[\bM_i\right]}{\partial \btheta_{i}(l)\partial \btheta_{i}(m)}\\
& =\frac{\partial^2 \vect\left[\bC_i\btheta_{i}\right]}{\partial \btheta_{i}(l)\partial \btheta_{i}(m)} =\bzero.
\end{split}
\end{equation}
Thus, the $(l,m)$-entry of $\bH_{\theta\theta,i}$ can be recast as
\begin{align}
& \bH_{\theta\theta,i}(l,m)=(K+1)\tr\left\{ \bX_i \frac{\partial \bM_i}{\partial \btheta_{i}(l)} 
\bX_i \frac{\partial \bM_i}{\partial \btheta_{i}(m)} \right\} \nonumber
\\
&-\tr\left\{
\bX_i (\bS_{\alpha}+\bS) \bX_i \frac{\partial \bM_i}{\partial \btheta_{i}(l)} \bX_i 
\frac{\partial \bM_i}{\partial \btheta_{i}(m)} 
\right\}\nonumber
\\
&
-\tr\left\{
\bX_i
(\bS_{\alpha}+\bS)
\bX_i
\frac{\partial \bM_i}{\partial \btheta_{i}(m)} 
\bX_i 
\frac{\partial \bM_i}{\partial \btheta_{i}(l)} 
\right\} \nonumber
\\
&=
(K+1)\tr\left\{ \bX_i  \frac{\partial \bM_i}{\partial \btheta_{i}(m)} 
\bF_i 
\bX_i \frac{\partial \bM_i}{\partial \btheta_{i}(l)} 
\right\} \nonumber
\\
&-\tr\left\{
\bX_i
(\bS_{\alpha}+\bS)
\bX_i
\frac{\partial \bM_i}{\partial \btheta_{i}(m)}
\bX_i 
\frac{\partial \bM_i}{\partial \btheta_{i}(l)} 
\right\}.
\end{align}
where
$$\bF_i = \left[ \bI_N - \bX_i \frac{(\bS_{\alpha}+\bS)}{(K+1)} \right]\, .$$
The above expression can be further simplified exploiting \eqref{eqn:traceProduct}. More precisely, the first term becomes
\begin{align}
& (K+1)\tr\left\{
\bX_i
\frac{\partial \bM_i}{\partial \btheta_{i}(m)}
\bF_i 
\bX_i 
\frac{\partial \bM_i}{\partial \btheta_{i}(l)} 
\right\} \nonumber
\\
&=
\begin{cases}
\ds (K+1) \left[\vect\left( \frac{\partial \bM_i}{\partial \btheta_{i}(l)} \right)  \right]^\dag
\left[\bX_i^* \bar{\bF}_i \otimes \bX_i\right]
\\
\times \ds \vect\left( \frac{\partial \bM_i}{\partial \btheta_{i}(m)} \right),
\\
\vspace{-4mm}
\\
\ds (K+1) \left[\vect\left( \frac{\partial \bM_i}{\partial \btheta_{i}(l)} \right)  \right]^T
\left[\bX_i \tilde{\bF}_i \otimes \bX_i\right] 
\\
\times \ds \vect\left( \frac{\partial \bM_i}{\partial \btheta_{i}(m)} \right),
\end{cases}\nonumber
\\
&=
\begin{cases}
\ds (K+1) \left[\bC_i \boe_{l,i}  \right]^\dag
\left[\bX_i^* \bar{\bF}_i \otimes \bX_i \right] \bC_i\boe_{m,i},
\\
\mbox{if } \bM_i \mbox{ is Hermitian},
\\
\ds (K+1) \left[\bC_i \boe_{l,i} \right]^T
[\bX_i \tilde{\bF}_i \otimes \bX_i] \bC_i \boe_{m,i},
\\
\mbox{if } \bM_i \mbox{ is symmetric}.
\end{cases}
\end{align}
where
$$\bar{\bF}_i = \left( \bI_N - \frac{((\bS_{\alpha}+\bS)\bX_i)^*}{K+1} \right) \, ,$$
$$\tilde{\bF}_i = \left( \bI_N - \frac{(\bS_{\alpha}+\bS)^*\bX_i}{K+1} \right) \, .$$
Using the same line of reasoning, it is possible to recast the last term as follows
\begin{align}
&\tr\left\{
\bX_i (\bS_{\alpha}+\bS) \bX_i 
\frac{\partial \bM_i}{\partial \btheta_{i}(m)} 
\bX_i 
\frac{\partial \bM_i}{\partial \btheta_{i}(l)} 
\right\}\nonumber
\\
&=
\begin{cases}
\ds \left[\vect\left( \frac{\partial \bM_i}{\partial \btheta_{i}(l)} \right)  \right]^\dag
[\bX_i^* \otimes \bX_i (\bS_{\alpha}+\bS) \bX_i ]\\ 
\times \vect\left( \frac{\partial \bM_i}{\partial \btheta_{i}(m)} \right)
\\
\vspace{-4mm}
\\
\ds \left[\vect\left( \frac{\partial \bM_i}{\partial \btheta_{i}(l)} \right)  \right]^T
[\bX_i \otimes \bX_i(\bS_{\alpha}+\bS)\bX_i  ]\\
\times \vect\left( \frac{\partial \bM_i}{\partial \btheta_{i}(m)} \right)
\end{cases}\nonumber
\\
&=
\begin{cases}
\ds \left[\bC_i \boe_{l,i}  \right]^\dag
[( \bX_i )^* \otimes  \bX_i  (\bS_{\alpha}+\bS)  \bX_i  ] \bC_i\boe_{m,i},
\\
\mbox{if } \bM_i \mbox{ is Hermitian},
\\
\ds \left[\bC_i \boe_{l,i}  \right]^T
[ \bX_i \otimes \bX_i  (\bS_{\alpha}+\bS)  \bX_i  ] \bC_i \boe_{m,i},
\\
\mbox{if } \bM_i \mbox{ is symmetric}.
\end{cases}
\end{align}
Summarizing, if $\bM_i$ is Hermitian $\bH_{\theta\theta,i}$ can be written as
\begin{align}
\bH_{\theta\theta,i}&=
(K+1) \bC_i^\dag
[ \bX_i^*  \bar{\bF}_i \otimes \bX_i ] \bC_i \nonumber
\\
&-\bC_i^\dag
[( \bX_i )^* \otimes  \bX_i(\bS_{\alpha}+\bS)\bX_i  ] \bC_i,
\end{align}
whereas if $\bM_i$ is symmetric we have that
\begin{align}
\bH_{\theta\theta,i}&=
(K+1) \bC_i^T
[ \bX_i \tilde{\bF}_i \otimes \bX_i ] \bC_i \nonumber
\\
&-\bC_i^T [ \bX_i \otimes  \bX_i(\bS_{\alpha}+\bS)\bX_i  ] \bC_i.
\end{align}
Next, consider $\bH_{\alpha\alpha,i}$ and observe that the gradient of \eqref{eqn:gradient_zAalpha} with respect to $\balpha^T$ is
\begin{equation}
\begin{split}
& \frac{\partial}{\partial \balpha^T}\left[ \frac{\partial s(\bp_i,H_i;\bz,\bZ)}{\partial \balpha} \right]
=\bH_{\alpha\alpha,i}\\
& =-2\left[
\begin{array}{cc}
\bv^\dag \bX_i\bv & \bzero
\\
\bzero & \bv^\dag \bX_i\bv
\end{array}
\right]
=-2 \bv^\dag \bX_i\bv \bI_2.
\end{split}
\end{equation}
As a final step towards the evaluation of $\bH_i$, we derive the expression for $\bH_{\alpha\theta,i}$. More precisely, 
exploiting previous results we get 
\begin{equation} \label{eqn:HalphaTheta}
\begin{split}
\frac{\partial }{\partial \btheta_i^T(l)} 
& \left[ \frac{\partial s(\bp_i,H_i;\bz,\bZ)}{\partial \balpha} \right]\\
& =\left[
\begin{array}{c}
\ds 2\alpha_{re} \tr\left[ \bA_1 \right]
-2 \Re\left\{\tr\left[ \bA_2 \right]\right\}
\\
\ds 2\alpha_{im} \tr\left[ \bA_1 \right]
+ 2 \Im\left\{\tr\left[ \bA_2	 \right]\right\}
\end{array}
\right],
\end{split}
\end{equation}
where $\bA_1=\bX_i \bv\bv^\dag \bX_i \frac{\partial \bX_i}{\partial \btheta_i(l)}$ and 
$\bA_2=\bX_i \bv\bz^\dag \bX_i \frac{\partial \bX_i}{\partial \btheta_i(l)}$.
Now, assume that the ICM is Hermitian; then using   \eqref{eqn:vecC}, \eqref{eqn:traceProduct}, and \eqref{eqn:vedDerC} the above equation can be recast as
\begin{equation}
\begin{split}
& \frac{\partial }{\partial \btheta_i^T(l)} \left[ \frac{\partial s(\bp_i,H_i;\bz,\bZ)}{\partial \balpha} \right]\\
& =\left[
\begin{array}{c}
\ds 2\alpha_{re} (\bC_i\boe_{l,i})^\dag \bar{\bPhi}_i
-2 \Re\left\{ (\bC_i\boe_{l,i})^\dag \tilde{\bPhi}_i \right\}
\\
\ds 2\alpha_{im} (\bC_i\boe_{l,i})^\dag \bar{\bPhi}_i
+2 \Im\left\{ (\bC_i\boe_{l,i})^\dag \tilde{\bPhi}_i \right\}
\end{array}
\right],
\end{split}
\end{equation}
where $\bar{\bPhi}_i = (\bX_i^*\otimes \bX_i) \vect[\bv\bv^\dag]$, and $\tilde{\bPhi}_i = (\bX_i^* \otimes \bX_i) \vect[\bv\bz^\dag]$. As a consequence,
\begin{equation}
\begin{split}
& \frac{\partial }{\partial \btheta_i^T} \left[ \frac{\partial s(\bp_i,H_i;\bz,\bZ)}{\partial \balpha} \right]\\
& =\left[
\begin{array}{c}
\left\{\ds 2\alpha_{re} \bC_i^\dag \bar{\bPhi}_i
-2 \Re\left\{ \bC_i^\dag \tilde{\bPhi}_i \right\}\right\}^T
\\
\left\{\ds +2\alpha_{im} \bC_i^\dag \bar{\bPhi}_i
+2 \Im\left\{ \bC_i^\dag \tilde{\bPhi}_i \right\}\right\}^T
\end{array}
\right].
\end{split}
\end{equation}
On the other hand, if the ICM is symmetric, then \eqref{eqn:HalphaTheta} becomes
\begin{equation}
\begin{split}
& \frac{\partial }{\partial \btheta_i^T(l)} \left[ \frac{\partial s(\bp_i,H_i;\bz,\bZ)}{\partial \balpha} \right]\\
& =\left[
\begin{array}{c}
\ds 2\alpha_{re} (\bC_i\boe_{l,i})^T \bar{\bPsi}_i
-2 \Re\left\{ (\bC_i\boe_{l,i})^T \tilde{\bPsi}_i \right\}
\\
\ds 2\alpha_{im} (\bC_i\boe_{l,i})^T \bar{\bPsi}_i
+2 \Im\left\{ (\bC_i\boe_{l,i})^T \tilde{\bPsi}_i \right\}
\end{array}
\right],
\end{split}
\end{equation}
where $\bar{\bPsi}_i = (\bX_i\otimes \bX_i) \vect[\bv\bv^\dag]$, and $\tilde{\bPsi}_i = (\bX_i\otimes \bX_i) \vect[\bv\bz^\dag]$. As a consequence,
\begin{equation}
\begin{split}
& \frac{\partial }{\partial \btheta_i^T} \left[ \frac{\partial s(\bp_i,H_i;\bz,\bZ)}{\partial \balpha} \right]\\
&=\left[
\begin{array}{c}
\left\{\ds 2\alpha_{re} \bC_i^T \bar{\bPsi}_i
-2 \Re\left\{ \bC_i^T \tilde{\bPsi}_i \right\}\right\}^T
\\
\left\{\ds 2\alpha_{im} \bC_i^T \bar{\bPsi}_i
+2 \Im\left\{ \bC_i^T \tilde{\bPsi}_i \right\}\right\}^T
\end{array}
\right].
\end{split}
\end{equation}


\bibliographystyle{IEEEtran}
\bibliography{MOS_paper}

\begin{table}[H] \centering
\begin{tabular}{|c|c|c|}
\hline
Parameter    & Case 1 ($L=1$) & Case 2 ($L=2$)\\
\hline
$N$   & 13 & 13\\
\hline
$\sigma_d$   & 0.15 & 0.15\\
\hline
$\rho_1$     & 0.85 & 0.85\\
\hline
$f_1$        & 0.285 & 0.285\\
\hline
$CNR_1$ [dB] & 30 & 20\\
\hline
$\rho_2$     & - & 0.93\\
\hline
$f_2$        & - & 0.05\\
\hline
$CNR_2$ [dB] & -  & 30\\
\hline
\end{tabular}
\caption{Parameters setting.}
\label{tab:case_study}
\end{table}
\end{document}

%% file: Defs2.tex
\newcommand{\bbm}{\begin{bmatrix}}
\newcommand{\ebm}{\end{bmatrix}}

\newcommand{\bit}{\begin{itemize}}
\newcommand{\eit}{\end{itemize}}

\newcommand{\ben}{\begin{enumerate}}
\newcommand{\een}{\end{enumerate}}

\newcommand{\bdesc}{\begin{description}}
\newcommand{\edesc}{\end{description}}

\newcommand{\bea}{\begin{array}}
\newcommand{\eea}{\end{array}}

\newcommand{\tr}{\mbox{\rm Tr}\, }

\newcommand{\beqa}{\begin{eqnarray}}
\newcommand{\eeqa}{\end{eqnarray}}

\newcommand{\ds}{\displaystyle}

\newcommand{\Comment}[1]{}

\def\H{{\mathcal H}}

\def\R{{\mathds R}}

\def\C{{\mathds C}}

\def\cC{\mbox{$\CMcal C$}}

\def\cN{\mbox{$\CMcal N$}}
\def\cO{\mbox{$\mathcal O$}}

\def\cU{\mbox{$\mathcal U$}}

\def\cbJ{\mbox{\mbox{\boldmath{$\CMcal J$}}}}
\def\cbI{\mbox{\mbox{\boldmath{$\CMcal I$}}}}

\newcommand{\be}{\begin{equation}}
\newcommand{\ee}{\end{equation}}

\newcommand{\bzero}{{\mbox{\boldmath $0$}}}

\newcommand{\boa}{{\mbox{\boldmath $a$}}}

\newcommand{\bc}{{\mbox{\boldmath $c$}}}

\newcommand{\boe}{{\mbox{\boldmath $e$}}}

\newcommand{\bm}{{\mbox{\boldmath $m$}}}
\newcommand{\bp}{\mbox{\boldmath $p$}}

\newcommand{\bv}{{\mbox{\boldmath $v$}}}
\newcommand{\bx}{{\mbox{\boldmath $x$}}}

\newcommand{\bz}{{\mbox{\boldmath $z$}}}

\newcommand{\bA}{{\mbox{\boldmath $A$}}}
\newcommand{\bB}{{\mbox{\boldmath $B$}}}
\newcommand{\bC}{{\mbox{\boldmath $C$}}}
\newcommand{\bD}{{\mbox{\boldmath $D$}}}

\newcommand{\bF}{{\mbox{\boldmath $F$}}}

\newcommand{\bH}{{\mbox{\boldmath $H$}}}
\newcommand{\bI}{{\mbox{\boldmath $I$}}}
\newcommand{\bJ}{{\mbox{\boldmath $J$}}}

\newcommand{\bM}{{\mbox{\boldmath $M$}}}

\newcommand{\bR}{{\mbox{\boldmath $R$}}}
\newcommand{\bS}{{\mbox{\boldmath $S$}}}

\newcommand{\bV}{{\mbox{\boldmath $V$}}}
\newcommand{\bX}{{\mbox{\boldmath $X$}}}

\newcommand{\bW}{{\mbox{\boldmath $W$}}}
\newcommand{\bZ}{{\mbox{\boldmath $Z$}}}

\newcommand{\diag}{\mbox{\boldmath\bf diag}\, }

\newcommand{\vect}{\mbox{\bf vec}\, }

\newcommand{\balpha}{{\mbox{\boldmath $\alpha$}}}

\newcommand{\btheta}{{\mbox{\boldmath $\theta$}}}

\newcommand{\bPhi}{{\mbox{\boldmath $\varPhi$}}}
\newcommand{\bPsi}{{\mbox{\boldmath $\Psi$}}}

